\documentclass[a4paper, 12pt] {article}
\usepackage[cp1251]{inputenc}
\usepackage[russian, english]{babel}
\usepackage{amsfonts}
\usepackage{mathtext}
\usepackage{cite}

\usepackage{graphicx}

\usepackage{verbatim}
\usepackage{amsmath}
\usepackage{amsthm}
\usepackage{amssymb}
\usepackage{delarray}
\usepackage{mathrsfs}
\usepackage{xcolor}

\begin{document}

\thispagestyle{empty}

\begin{center}
{\bf\large AN INVERSE SPECTRAL PROBLEM FOR SECOND-ORDER FUNCTIONAL-DIFFERENTIAL PENCILS WITH TWO DELAYS}
\end{center}

\begin{center}
{\bf\large S.A. Buterin\footnote{Department of Mathematics, Saratov State University, Russia, {\it email:
buterinsa@info.sgu.ru}}, M.A. Malyugina\footnote{Department of Mathematics, Saratov State University, Russia,
{\it email: margarita.malyugina@tele2.ru}} and C.-T. Shieh\footnote{Department of Mathematics, Tamkang
University, Taiwan, {\it email: ctshieh@mail.tku.edu.tw}}}
\end{center}

{\bf Abstract.} We consider a second order functional-differential pencil with two constant delays of the
argument and study the inverse problem of recovering its coefficients from the spectra of two boundary value
problems with one common boundary condition. The uniqueness theorem is proved and a constructive procedure
for solving this inverse problem along with necessary and sufficient conditions for its solvability is
obtained. Moreover, we give a survey on the contemporary state of the inverse spectral theory for operators
with delay. The pencil under consideration generalizes Sturm--Liouville-type operators with delay, which
allows us to illustrate essential results in this direction, including recently solved open questions.


\\

{\large\bf 1. Introduction and main results}
\\

Recently, there appeared considerable interest in inverse problems of spectral analysis for
Sturm--Liouville-type operators with constant delay:
\begin{equation}\label{0}
\ell y\equiv-y''(x)+q(x)y(x-a)=\lambda y(x),\quad 0<x<\pi,
\end{equation}
under two-point boundary conditions, see \cite{Pik91, FrYur12, VladPik16, ButYur19, ButPikYur17, Ign18,
BondYur18-1, BondYur18-2, VPV19, Yur19, DV19, Yur19-2, Yur20, SS19, DurBut, Yang14, WShM19}, which are often
adequate for modelling various real-world processes frequently possessing a nonlocal nature. Here $q(x)$ is a
complex-valued function in $L_2(a,\pi)$ vanishing on $(0,a).$ In particular, it is well known that
specification of the spectra $\{\lambda_{n,j}\},$ $j=0,1,$ of two boundary value problems for the
functional-differential equation~(\ref{0}) with one common boundary condition in zero, say,
\begin{equation}\label{2}
y(0)=y^{(j)}(\pi)=0
\end{equation}
uniquely determines the potential $q(x)$ as soon as $a\in[\pi/2,\pi).$ Moreover, the corresponding inverse
problem is overdetermined. Thus, in \cite{ButYur19}, conditions on an arbitrary increasing sequence of
natural numbers $\{n_k\}_{k\ge1}$ were obtained that are necessary and sufficient for the unique
determination of $q(x)$ by specifying the corresponding subspectra $\{\lambda_{n_k,0}\}$ and
$\{\lambda_{n_k,1}\}.$

For a long time, it was an open question whether the uniqueness result for $a\in[\pi/2,\pi)$ would remain
true also for $a\in(0,\pi/2).$ The positive answer for $a\in[2\pi/5,\pi/2)$ was given independently in
\cite{BondYur18-1}, and in \cite{VPV19} for the case of the Robin boundary condition in zero:
\begin{equation}\label{rob}
y'(0)-hy(0)=y^{(j)}(\pi)=0.
\end{equation}

For $a\in[\pi/3,2\pi/5),$ the authors of \cite{BondYur18-2} have shown that the spectra of both problems
consisting of (\ref{0}) and (\ref{2}) uniquely determine the potential $q(x)$ on $(a,3a/2)\cup(\pi-a/2,\pi).$
But the strongest uniqueness result under these settings was obtained in \cite{DV19}, where it was proved
that $q(x)$ is uniquely determined on the set $(a,3a/2)\cup(\pi-a,2a)\cup(\pi-a/2,\pi).$ Moreover, the
authors of \cite{DV19} proved that, for the complete determination of $q(x),$ it is sufficient to
additionally specify it on $(3a/2,\pi/2+a/4)$ as well as its mean value on $(\pi/2+a/4,\pi-a).$

Meanwhile, the recent paper \cite{DurBut} gave a negative answer to the open question formulated above for
boundary conditions~(\ref{2}) as soon as $a\in[\pi/3,2\pi/5)$ by constructing a one-parametric infinite
family $B$ of different iso-bispectral potentials~$q(x),$ i.e. of those for which both boundary value
problems possess one and the same pair of spectra. It is also interesting that potentials in $B$ differ on
the set $(3a/2,\pi-a)\cup(2a,\pi-a/2)$ by an arbitrary multiplicative complex constant. Thus, the uniqueness
subdomain for $q(x)$ established in \cite{DV19} is unimprovable.

This appeared quite unexpected taking into account that the paper \cite{Yur20} announced that specification
of the spectra of both boundary value problems consisting of~(\ref{0}) and~(\ref{rob}) uniquely determines
the potential $q(x)$ also for $a\in[\pi/3,2\pi/5),$ see also \cite{Yur19}. Even though, in \cite{Yur20},
Robin boundary conditions were imposed also in the point $\pi,$ they can be easily reduced to~(\ref{rob}).

Moreover, for any fixed $a\in(0,\pi),$ in the papers \cite{FrYur12} and \cite{ButPikYur17} for the cases of
boundary conditions (\ref{2}) and (\ref{rob}), respectively, it was established that if the spectra coincide
with the ones of the corresponding problems with the zero potential, then $q(x)$ is zero too.

Among numerous studies devoted to inverse spectral problems for functional-differential operators with delay,
to the best of our knowledge, \cite{ButYur19} remains a sole work dealing with necessary and sufficient
conditions of solvability for inverse problems of this class. In particular, from the results in
\cite{ButYur19} it follows that, unlike the classical case $a=0,$ the spectra $\{\lambda_{n,0}\}$ and
$\{\lambda_{n,1}\}$ may have any finite number of common eigenvalues. In Appendix~A, we supplement the work
\cite{ButYur19}, by applying the solvability result of the present paper to (\ref{0}) for $a\in[2\pi/5,\pi).$

There are also works devoted to inverse problems for operators with several delays:
\begin{equation}\label{0.1}
\ell_m y\equiv-y''(x)+\sum_{\nu=1}^m q_\nu(x)y(x-a_\nu)=\lambda y(x),\quad 0<x<\pi, \quad m>1,
\end{equation}
where $q_\nu(x)=0$ on $(0,a_\nu)$ for $\nu=\overline{1,n},$ see \cite{PVP16, Sh19, Sh20, VPV20, VPVC20}.
However, all potentials $q_\nu(x)$ cannot be completely determined simultaneously even by specifying
arbitrarily many different spectra. This becomes especially obvious when all $a_\nu$ are equal, but even for
different $a_\nu$ there remain subintervals, where the functions $q_\nu(x)$ cannot be distinguished (see
Appendix~B).

An attempt to generalize the results of papers \cite{ButPikYur17, FrYur12} to the operator $\ell_2$ was made
in \cite{Sh19}, where it was claimed that if the spectra of two boundary value problems for one and the same
equation~(\ref{0.1}) for $m=2$ along with the boundary conditions (\ref{2}) coincide with the analogous two
spectra corresponding to the zero potentials, then $q_1(x)$ and $q_2(x)$ are zeros too. Thereafter, the
generalization of this assertion to arbitrary $m>1$ was announced in \cite{Sh20}. But, unfortunately, both
papers \cite{Sh19, Sh20} contain a serious mistake, and these two assertions {\it cannot} be true even for
different delays (see the counterexample in Appendix~B).

Nevertheless, one can completely recover all potentials $q_\nu(x)$ by using the spectra of $2m$ boundary
value problems for~$m$ different equations specially composed from (\ref{0.1}). For example, in \cite{VPV20,
VPVC20} for $m=2$ and $a_1,a_2\in[\pi/2,\pi),$ it was shown that, for recovering $q_1(x)$ and $q_2(x),$ it is
sufficient to specify the spectra of four boundary value problems for two equations:
$$
-y''(x)+q_1(x)y(x-a_1) +(-1)^\nu q_2(x)y(x-a_2)=\lambda y(x),\quad 0<x<\pi, \quad \nu=1,2.
$$

In the present paper, we show, in particular, that there is no need to use two different equations if one of
delayed terms depends on the spectral parameter. Specifically, we consider the functional-differential
equation with nonlinear dependence on the spectral parameter $\rho:$
\begin{equation}\label{1}
y''(x)+\rho^2y(x)=q_0(x)y(x-a_0)+2\rho q_1(x)y(x-a_1),\quad 0<x<\pi,
\end{equation}
which generalizes equation (\ref{0}). We assume that $a_0\in[\pi/3,\pi),$ $a_1\in[\pi/2,\pi)$ and
$a_0+a_1\ge\pi.$ The case $a_0+a_1<\pi$ requires a separate investigation (see Section~3 for details). For
$\nu=0,1,$ let $q_\nu(x)$ be a complex-valued function in $W_2^\nu[a_\nu,\pi],$ $q_\nu(x)=0$ on $(0,a_\nu)$
and
\begin{equation}\label{1.1}
\int\limits_{a_1}^\pi q_1(x)\,dx=0.
\end{equation}
For $j=0,1,$ we denote by $\{\rho_{n,j}\}$ the spectrum of the boundary value problem ${\cal L}_j:={\cal
L}_j(q_0,q_1)$ that consists of equation~(\ref{1}) along with the boundary conditions (\ref{2}). We also use
one and the same symbol $\{\varkappa_n\}$ for denoting {\it different} sequences in~$l_2,$ and put ${\mathbb
Z}_0:={\mathbb Z}\setminus\{0\}$ and ${\mathbb Z}_1:={\mathbb Z}.$ Under~our assumptions, we prove first the
following theorem giving asymptotics of the spectra.

\medskip
{\bf Theorem 1.} {\it For $j=0,1$ and $n\in{\mathbb Z}_j,$ the following asymptotics holds:
\begin{equation}\label{2.1}
\rho_{n,j}=\rho_{n,j}^0 +\frac{\omega}{\pi n}\cos\rho_{n,j}^0a_0 +\frac{\alpha_j}{\pi n}\sin\rho_{n,j}^0a_1
+\frac{\varkappa_n}n, \quad \alpha_0,\alpha_1,\omega\in{\mathbb C},
\end{equation}
where $\rho_{n,j}^0=n-j/2.$ Moreover,
\begin{equation}\label{2.2}
\alpha_j=\alpha+(-1)^j\beta, \quad \alpha=\frac{q_1(a_1)}2, \quad \beta=\frac{q_1(\pi)}2,
\quad\omega=\frac12\int\limits_{a_0}^\pi q_0(x)\,dx.
\end{equation}}

For $\alpha=0,$ this theorem was announced in the conference papers \cite{BM19,BMS19}, which became the first
works dealing with inverse problems for functional-differential pencils in any form.

Consider the following inverse problem.

\medskip
{\bf Inverse Problem 1.} Given the spectra $\{\rho_{n,j}\}_{n\in{\mathbb Z}_j},\,j=0,1;$ find $q_0(x)$ and
$q_1(x).$

\medskip
We note that for $a_0=a_1=0$ this inverse problem was studied in \cite{ButYur12}. The following theorem gives
uniqueness of its solution under our present settings.

\medskip
{\bf Theorem 2. }{\it Let both spectra $\{\rho_{n,j}\}_{n\in{\mathbb Z}_j},\,j=0,1,$ be specified. Then the
function $q_0(x)$ is uniquely determined a.e. on the union of intervals $I_1:=(a_0,3a_0/2)\cup
(\pi-a_0,2a_0)\cup (\pi-a_0/2,\pi),$ while the function $q_1(x)$ is uniquely determined on the entire segment
$[a_1,\pi].$}

\medskip
In particular, both functions $q_0(x)$ and $q_1(x)$ are completely determined when $a_0\ge2\pi/5.$ However,
for $a_0\in[\pi/3,2\pi/5),$ the complete uniqueness does not take place. For illustrating this, one can use
the same one-parametric family of functions $B=\{q_{0,\gamma}(x)\}_{\gamma\in{\mathbb C}}$ constructed in
\cite{DurBut}, for which also the problems ${\cal L}_j(q_{0,\gamma},q_1),\,j=0,1,$ will possess one and the
same pair of spectra for all $\gamma\in{\mathbb C}.$ Since, as was already mentioned above, for different
values of $\gamma$ the functions $q_{0,\gamma}(x)$ differ precisely on the set $(a_0,\pi) \setminus{\rm
cl}(I_1)$ except their common zeros, the uniqueness subdomain~$I_1$ cannot be refined.

However, as in \cite{DV19}, the function $q_0(x)$ would be determined uniquely if some a priori information
on it were additionally specified. Namely, the following theorem holds.

\medskip
{\bf Theorem 3. }{\it Under the hypothesis of Theorem~2, the specification of the function $q_0(x)$ on the
subinterval $(3a_0/2,\pi/2+a_0/4)$ along with the value
$$
\omega_0:=\int\limits_{\frac\pi2+\frac{a_0}4}^{\pi-a_0} q_0(x)\,dx.
$$
determines it also on $(\pi/2+a_0/4,\pi-a_0)\cup(2a_0,\pi-a_0/2).$ So $q_0(x)$ is determined completely.}

\medskip
As in \cite{ButYur19}, one can show that, in the case $a_0\ge2\pi/5,$ Inverse Problem~1 is overdetermined,
and also describe subspectra, whose specification would uniquely determine the functions $q_0(x)$ and
$q_1(x).$ Some results of this type can be found in \cite{BM19,BMS19}. But here we restrict ourself to
dealing with the full spectra, which surprisingly does not prevent us from obtaining necessary and sufficient
conditions for the solvability of Inverse Problem~1 for $a_0\ge2\pi/5.$ Besides asymptotics~(\ref{2.1}),
these conditions include some restrictions on the growth of certain entire functions, which makes the inverse
problem consistent even in spite of its overdetermination. Specifically, the following theorem~holds.

\medskip
{\bf Theorem 4. }{\it Let $a_0\ge2\pi/5.$ Then for any sequences of complex numbers
$\{\rho_{n,0}\}_{|n|\in{\mathbb N}}$ and $\{\rho_{n,1}\}_{n\in{\mathbb Z}}$ to be the spectra of some
boundary value problems ${\cal L}_0(q_0,q_1)$ and ${\cal L}_1(q_0,q_1),$ respectively, it is necessary and
sufficient to satisfy the following two conditions:

(i) For $j=0,1,$ the sequence $\{\rho_{n,j}\}_{n\in{\mathbb Z}_j}$ has the form~(\ref{2.1});

(ii) For $j,\nu=0,1,$ the exponential type of the function $g_{j,\nu}(\rho)$ does not exceed $\pi-a_\nu,$
where
\begin{equation}\label{4.32}
g_{j,\nu}(\rho)=\theta_j(\rho)+(-1)^{j+\nu}\theta_j(-\rho),
\end{equation}
\begin{equation}\label{4.33}
\theta_0(\rho)=\rho^2\Delta_0(\rho) -\rho\sin\rho\pi +\omega\cos\rho(\pi-a_0) -\alpha_0 \sin\rho(\pi-a_1),
\end{equation}
\begin{equation}\label{4.34}
\theta_1(\rho)=\rho\Delta_1(\rho) -\rho\cos\rho\pi -\omega\sin\rho(\pi-a_0) -\alpha_1\cos\rho(\pi-a_1),
\end{equation}
while the functions $\Delta_j(\rho)$ are determined by the formula
\begin{equation}\label{4.13}
\Delta_j(\rho)=\pi^{1-j}\prod_{n\in{\mathbb Z}_j}\frac{\rho_{n,j}-\rho}{\rho_{n,j}^0}
\exp\Big(\frac\rho{\rho_{n,j}^0}\Big), \quad j=0,1.
\end{equation}}

For $j=0,1,$ the function $\Delta_j(\rho)$ determined by (\ref{4.13}) is the characteristic function of the
problem ${\cal L}_j$ (see the next section). The proof of Theorem~4 is constructive and gives an algorithm
for solving Inverse Problem~1 (Algorithm~1 in Section~6). For proving Theorem~4, we obtain and study a
transformation operator associated with equation (\ref{0}), which allows reducing the inverse problem to the
so-called main vectorial integral equation. For applications of the transformation operator approach to other
classes of nonlocal operators, see the survey~\cite{But20}.

The paper is organized as follows. In the next section, we construct a transformation operator for the
sine-type solution of equation (\ref{1}), and study the characteristic functions of the problems ${\cal
L}_j.$ Therein, we also give the proof of Theorem~1. In Section~3, we derive and study the main equation of
Inverse Problem~1. In Section~4, we prove the uniqueness theorems (Theorems~2 and~3). In Section~5, we obtain
important representations for the functions determined by (\ref{4.13}) with arbitrary complex zeros of the
form (\ref{2.1}). In Section~6, we prove Theorem~4 and obtain an algorithm for solving the inverse problem.
In Appendix~A, we provide an analog of Theorem~4 for the Sturm--Liouville-type operator $\ell$ with the delay
$a\in[2\pi/5,\pi),$ i.e. when $q_1(x)\equiv0.$ In Appendix~B, we give a counterexample, showing that
specification of any spectra does not uniquely determine all potentials in equation (\ref{0.1}) even for
$m=2$ and $a_1\ne a_2.$
\\

{\large\bf 2. Transformation operator and characteristic functions}
\\

Let $y=S(x,\rho)$ be the sine-type solution of equation (\ref{1}), i.e. the solution satisfying the initial
conditions $S(0,\rho)=0$ and $S'(0,\rho)=1.$ By virtue of its uniqueness, eigenvalues of the problem ${\cal
L}_j,$ $j=0,1,$ coincide with zeros of the entire function
\begin{equation}\label{3}
\Delta_j(\rho):=S^{(j)}(\pi,\rho),
\end{equation}
which is called {\it characteristic function} of ${\cal L}_j.$ We introduce the designations
$$
Q_\nu(x):=\int\limits_{a_\nu}^x q_\nu(t)\,dt,\;\; \nu=0,1, \quad c_0(x):=\cos x, \quad c_1(x):=\sin x,
$$
where the latter two ones are aimed to be used only occasionally for convenience.

The following lemma gives the transformation operator that connects the solutions $\rho^{-1}\sin\rho x$ and
$\cos\rho x$ of the simplest equation (\ref{1}), possessing zero coefficients, with the solution $S(x,\rho).$

\medskip
{\bf Lemma 1. }{\it The following representation holds:
\begin{equation}\label{4}
S(x,\rho)=\frac{\sin\rho x}{\rho} -Q_1(x) \frac{\cos\rho(x-a_1)}{\rho} +\sum_{\nu=0}^1\int\limits_{a_\nu}^x
K_\nu(x,t)\frac{c_{1-\nu}(\rho(x-t))}\rho\,dt, \quad 0\le x\le\pi,
\end{equation}
where $K_\nu(x,t)=0$ in the exterior of the triangle $a_\nu\le t\le x \le\pi$ as $\nu=0,1,$ and
\begin{equation}\label{4.1}
K_1(x,t)=\frac12\Big(q_1\Big(x-\frac{t-a_1}2\Big) +q_1\Big(\frac{t+a_1}2\Big)\Big), \quad a_1\le t\le
x\le\pi,
\end{equation}
while the kernel $K_0(x,t)$ satisfies the following integral equation:
\begin{equation}\label{4.2}
K_0(x,t)=\frac12\int\limits_\frac{t+a_0}2^{x-\frac{t-a_0}2}q_0(\tau)\,d\tau +\frac{A(x,t)}2, \quad a_0\le
t\le x \le\pi,
\end{equation}
where
\begin{equation}\label{4.3}
A(x,t)=\left\{\begin{array}{l} 0, \qquad a_0\le t \le\min\{x,2a_0\}, \;\; x\le\pi,\\[3mm]
\displaystyle \int\limits_t^xq_0(\tau)\,d\tau\int\limits_{a_0}^{t-a_0} K_0(\tau-a_0,\eta)\,d\eta - \!\!\!
\int\limits_{x-\frac{t}2+a_0}^xq_0(\tau)\,d\tau
\!\!\!\int\limits_{a_0}^{2(\tau-x)+t-a_0} \!\!\! K_0(\tau-a_0,\eta)\,d\eta\\[3mm]
\qquad\qquad \displaystyle +\int\limits_{\frac{t}2+a_0}^tq_0(\tau)\,d\tau\int\limits_{a_0}^{2\tau-t-a_0}
K_0(\tau-a_0,\eta)\,d\eta, \qquad 2a_0< t\le x\le\pi.
\end{array}\right.
\end{equation}
In particular, the following relations hold:
\begin{equation}\label{4.4}
K_0(x,x)=0, \;\; K_0(x,a_0)=\frac{Q_0(x)}2, \;\; K_1(x,a_1)=\frac{q_1(x)}2+\alpha, \;\;
A(\pi,2a_0)=A(\pi,\pi)=0.
\end{equation}
}

{\bf Remark 1.} Since $\pi\le3a_0,$ the integration variable $\eta$ in (\ref{4.3}) never exceeds $2a_0.$
Thus, the function $K_0(\tau-a_0,\eta)$ under the integrals in (\ref{4.3}) is determined by the formula
\begin{equation}\label{4.3.1}
K_0(\tau-a_0,\eta)=\frac12\int\limits_\frac{\eta+a_0}2^{\tau-\frac{\eta+a_0}2}q_0(\zeta)\,d\zeta, \quad
a_0\le \eta \le\min\{\tau-a_0,2a_0\}, \;\; \tau\le\pi.
\end{equation}

\medskip
{\it Proof of Lemma~1.} Clearly, the function $S(x,\rho)$ obeys the following integral equation:
\begin{equation}\label{4.5}
S(x,\rho)=\frac{\sin\rho x}{\rho} +\sum_{\nu=0}^1 (2\rho)^\nu\int\limits_{a_\nu}^x \frac{\sin\rho(x-t)}\rho
q_\nu(t)S(t-a_\nu,\rho)\,dt, \quad 0\le x\le\pi.
\end{equation}
Since $a_\nu+a_1\ge\pi$ for $\nu=0,1,$ substituting (\ref{4}) into (\ref{4.5}), we arrive at the relation
\begin{equation}\label{4.6}
\int\limits_{a_0}^x K_0(x,t)\frac{\sin\rho(x-t)}\rho\,dt +\int\limits_{a_1}^x
K_1(x,t)\frac{\cos\rho(x-t)}\rho\,dt-Q_1(x) \frac{\cos\rho(x-a_1)}{\rho} =\sum_{j=0}^2 {\cal A}_j(x,\rho),
\end{equation}
where
$$
{\cal A}_\nu(x,\rho)=2^\nu\rho^{\nu-1}\int\limits_{a_\nu}^x \sin\rho(x-t) q_\nu(t)\,dt
\int\limits_0^{t-a_\nu} \cos\rho\tau\,d\tau, \quad \nu=0,1,
$$
$$
{\cal A}_2(x,\rho)=\int\limits_{2a_0}^x \frac{\sin\rho(x-t)}\rho q_0(t)\,dt\int\limits_{a_0}^{t-a_0}
K_0(t-a_0,\tau)\,d\tau\int\limits_0^{t-\tau-a_0} \cos\rho\xi\,d\xi.
$$
Using the formula $2\sin\rho(x-t)\cos\rho\tau =\sin\rho(x-t+\tau)+ \sin\rho(x-t-\tau)$ and changing the
integration variables along with its order, we get
\begin{equation}\label{4.7}
{\cal A}_\nu(x,\rho)=(2\rho)^{\nu-1}\int\limits_{a_\nu}^x \sin\rho(x-t)\,dt
\int\limits_\frac{t+a_\nu}2^{x-\frac{t-a_\nu}2} q_\nu(\tau)\,d\tau, \quad \nu=0,1,
\end{equation}
\begin{equation}\label{4.8}
{\cal A}_2(x,\rho)= \int\limits_{2a_0}^x A(x,t)\frac{\sin\rho(x-t)}{2\rho}\,dt,
\end{equation}
where $A(x,t)$ is determined by formula (\ref{4.3}). Moreover, integration by parts gives
\begin{equation}\label{4.9}
{\cal A}_1(x,\rho)= -Q_1(x) \frac{\cos\rho(x-a_1)}{\rho} +\int\limits_{a_1}^x
\Big(q_1\Big(x-\frac{t-a_1}2\Big)+q_1\Big(\frac{t+a_1}2\Big)\Big)\frac{\cos\rho(x-t)}{2\rho}\,dt.
\end{equation}
Substituting (\ref{4.7}) for $\nu=0$ as well as (\ref{4.8}) and (\ref{4.9}) into (\ref{4.6}), we arrive at
(\ref{4})--(\ref{4.3}). Relations (\ref{4.4}) are obvious. $\hfill\Box$

\medskip
The next lemma gives fundamental representations for the characteristic functions.

\medskip
{\bf Lemma 2. }{\it The following representations hold:
\begin{equation}\label{4.11}
\Delta_0(\rho)=\frac{\sin\rho\pi}\rho-\omega\frac{\cos\rho(\pi-a_0)}{\rho^2}
+\alpha_0\frac{\sin\rho(\pi-a_1)}{\rho^2} +\sum_{\nu=0}^1\int\limits_0^{\pi-a_\nu} w_{0,\nu}(x)
\frac{c_\nu(\rho x)}{\rho^2}\,dx,
\end{equation}
\begin{equation}\label{4.12}
\Delta_1(\rho)=\cos\rho\pi+\omega\frac{\sin\rho(\pi-a_0)}\rho +\alpha_1\frac{\cos\rho(\pi-a_1)}\rho
+\sum_{\nu=0}^1\int\limits_0^{\pi-a_\nu} w_{1,\nu}(x) \frac{c_{1-\nu}(\rho x)}\rho\,dx,
\end{equation}
where $w_{j,\nu}(x)\in L_2(0,\pi-a_\nu),$ $j,\nu=0,1,$ and
\begin{equation}\label{4.12.1.2}
\int\limits_0^{\pi-a_0} w_{0,0}(x)\,dx=\omega, \;\; \int\limits_0^{\pi-a_1}
xw_{0,1}(x)\,dx=\alpha_0(a_1-\pi),\;\; \int\limits_0^{\pi-a_1} w_{1,1}(x)\,dx=-\alpha_1.
\end{equation}
Moreover,
\begin{equation}\label{4.10}
w_{0,\nu}(x) =(-1)^{\nu+1}K_{\nu,2}(\pi,\pi-x), \quad w_{1,\nu}(x) =P_\nu(\pi,\pi-x), \quad \nu=0,1,
\end{equation}
where
\begin{equation}\label{4.10.1}
K_{\nu,1}(x,t):=\frac{\partial}{\partial x}K_\nu(x,t), \;\; K_{\nu,2}(x,t):=\frac{\partial}{\partial
t}K_\nu(x,t), \;\; P_\nu(x,t):=K_{\nu,1}(x,t)+K_{\nu,2}(x,t).
\end{equation}
}

\medskip
{\it Proof.} Integrating by parts in (\ref{4}) with account of (\ref{4.4}), and recalling (\ref{1.1}),
(\ref{2.2}) along with (\ref{3}), we get (\ref{4.11}). Differentiating (\ref{4}) with respect to $x$ and then
using integration by parts, we obtain~(\ref{4.12}). Finally, although relations (\ref{4.12.1.2}) can be
established by direct calculations, we accept them just as a simple corollary from entireness of the
functions $\Delta_0(\lambda)$ and $\Delta_1(\lambda).$ $\hfill\Box$

\medskip
Now we are in position to give the proof of Theorem~1.

\medskip
{\it Proof of Theorem~1.} By the standard approach involving Rouch\'e's theorem (see, e.g., \cite{FY01}),
using representation (\ref{4.11}), one can show that the function $\Delta_0(\lambda)$ has infinitely many
zeros of the form $\rho_{n,0}=n+\varepsilon_{n,0},$ where $|n|\in{\mathbb N},$ while $\varepsilon_{n,0}\to0$
as $|n|\to\infty.$ Substituting this representation into (\ref{4.11}), we obtain
\begin{equation}\label{4.10.2}
\sin\rho_{n,0}\pi=\frac{\omega}n\cos(n+\varepsilon_{n,0})(\pi-a_0)
-\frac{\alpha_0}n\sin(n+\varepsilon_{n,0})(\pi-a_1) +\frac{\varkappa_n}n.
\end{equation}
Since $\sin\rho_{n,0}\pi=\sin(n+\varepsilon_{n,0})\pi=(-1)^n\varepsilon_{n,0}\pi+O(\varepsilon_{n,0}^3)$ as
$|n|\to\infty,$ we refine $\varepsilon_{n,0}=O(n^{-1})$ for $|n|\to\infty.$ Hence, we have the asymptotic
formulae
$$
\sin\rho_{n,0}\pi=(-1)^n\varepsilon_{n,0}\pi+O\Big(\frac1{n^3}\Big),
$$
$$
\cos(n+\varepsilon_{n,0})(\pi-a_0)=(-1)^n\cos na_0+O\Big(\frac1n\Big),
$$
$$
\sin(n+\varepsilon_{n,0})(\pi-a_1)=(-1)^{n+1}\sin na_1+O\Big(\frac1n\Big)
$$
as soon as $|n|\to\infty.$ Substituting them into (\ref{4.10.2}), we arrive at
$$
\varepsilon_{n,0}=\frac{\omega}{\pi n}\cos na_0 +\frac{\alpha_0}{\pi n}\sin na_1 +\frac{\varkappa_n}n,
$$
which implies (\ref{2.1}) for $j=0.$

Analogously, applying Rouch\'e's theorem to representation (\ref{4.12}), we get
$\rho_{n,1}=\rho_{n,1}^0+\varepsilon_{n,1},$ where $n\in{\mathbb Z},$ while $\varepsilon_{n,1}\to0$ as
$|n|\to\infty.$ Substituting this into (\ref{4.12}), we obtain
\begin{equation}\label{4.10.3}
\cos\rho_{n,1}\pi=-\frac{\omega}n\sin(\rho_{n,1}^0+\varepsilon_{n,1})(\pi-a_0)
-\frac{\alpha_1}n\cos(\rho_{n,1}^0+\varepsilon_{n,1})(\pi-a_1) +\frac{\varkappa_n}n.
\end{equation}
Since $\cos\rho_{n,1}\pi=\cos(\rho_{n,1}^0+\varepsilon_{n,1})\pi=(-1)^n\varepsilon_{n,1}\pi
+O(\varepsilon_{n,1}^3)$ as $|n|\to\infty,$ we refine $\varepsilon_{n,1}=O(n^{-1})$ for $|n|\to\infty.$ Then
substituting the asymptotic formulae
$$
\sin(\rho_{n,1}^0+\varepsilon_{n,1})(\pi-a_0)=(-1)^{n+1}\cos\rho_{n,1}^0a_0 +O\Big(\frac1n\Big),
$$
$$
\cos(\rho_{n,1}^0+\varepsilon_{n,1})(\pi-a_1)=(-1)^{n+1}\sin\rho_{n,1}^0a_1 +O\Big(\frac1n\Big),
$$
for $|n|\to\infty$ into (\ref{4.10.3}), we arrive at
$$
\varepsilon_{n,1}=\frac{\omega}{\pi n}\cos \rho_{n,1}^0a_0 +\frac{\alpha_1}{\pi n}\sin \rho_{n,1}^0a_1
+\frac{\varkappa_n}n,
$$
which implies (\ref{2.1}) for $j=1.$  $\hfill\Box$

\medskip
Finally, we obtain formulae for recovering the characteristic functions from their zeros.

\medskip
{\bf Lemma 3. }{\it For any $a_0,a_1\in[0,2\pi],$ each function $\Delta_0(\rho)$ and $\Delta_1(\rho)$ of the
form described in (\ref{4.11})--(\ref{4.12.1.2}) is determined by its zeros uniquely. Moreover,
representation (\ref{4.13}) holds.}

\medskip
{\it Proof.} By virtue of Hadamard's factorization theorem (see, e.g., \cite{BFY14}), we get
\begin{equation}\label{4.14}
\Delta_j(\rho)=C_j\rho^{s_j}\exp(b_j\rho)\prod_{\rho_{n,j}\ne0}\Big(1-\frac\rho{\rho_{n,j}}\Big)
\exp\Big(\frac\rho{\rho_{n,j}}\Big), \quad j=0,1,
\end{equation}
where $C_j$ and $b_j$ are some constants, while $s_j$ is the multiplicity of the null zero $\rho_{n,j}=0.$ In
particular, we have
\begin{equation}\label{4.15}
\rho^{j-1}c_{1-j}(\rho\pi) =\pi^{1-j}\prod_{n\in{\mathbb Z}_j}\Big(1-\frac\rho{\rho_{n,j}^0} \Big)
\exp\Big(\frac\rho{\rho_{n,j}^0}\Big), \quad j=0,1.
\end{equation}
Dividing (\ref{4.14}) by (\ref{4.15}), we obtain
$$
\frac{\rho^{1-j}\Delta_j(\rho)}{c_{1-j}(\rho\pi)}=\frac{C_{j}}{\pi^{1-j}} \exp\Big(\Big(b_j
+\sum_{\rho_{n,j}\ne0}\Big(\frac1{\rho_{n,j}}-\frac1{\rho_{n,j}^0}\Big)-\sum_{\rho_{n,j}=0}\frac1{\rho_{n,j}^0}\Big)\rho\Big)
\qquad\qquad\qquad\qquad
$$
\begin{equation}\label{4.16}
\qquad\qquad\qquad\qquad \qquad\times\prod_{\rho_{n,j}=0}\Big(\frac1\rho-\frac1{\rho_{n,j}^0}\Big)^{-1}
\prod_{\rho_{n,j}\ne0}\frac{\rho_{n,j}^0}{\rho_{n,j}}
\prod_{\rho_{n,j}\ne0}\frac{\rho_{n,j}-\rho}{\rho_{n,j}^0-\rho}, \quad j=0,1.
\end{equation}
On the other hand, (\ref{4.11}) and (\ref{4.12}) imply $\rho^{1-j}(c_{1-j}(\rho\pi))^{-1} \Delta_j(\rho)\to1$
for $j=0,1$ as $\rho^2\to-\infty,$ which along with (\ref{4.16}) gives
$$
C_j=\pi^{1-j}(-1)^{s_j}\prod_{\rho_{n,j}=0}\frac1{\rho_{n,j}^0}
\prod_{\rho_{n,j}\ne0}\frac{\rho_{n,j}}{\rho_{n,j}^0}, \quad b_j=\sum_{\rho_{n,j}=0}\frac1{\rho_{n,j}^0}
+\sum_{\rho_{n,j}\ne0}\Big(\frac1{\rho_{n,j}^0} -\frac1{\rho_{n,j}}\Big).
$$
Substituting this into (\ref{4.14}), we arrive at (\ref{4.13}). $\hfill\Box$
\\

{\large\bf 3. Main equation of the inverse problem}
\\

The relations in (\ref{4.10}) can be considered as a system of equations with respect to the functions
$q_0(x)$ and $p(x):=q_1'(x),$ which we refer to as {\it main (vectorial) equation} of Inverse Problems~1. By
virtue of our standing assumption $a_0+a_1\ge\pi,$ for each $\nu\in\{0,1\},$ the functions $K_{\nu,2}(x,t)$
and $P_\nu(x,t)$ depend only on $q_\nu(x),$ while for $\nu=1$ they depend even only on $p(x).$ Hence, the
main equation can be splitted into two independent subsystems for $\nu=0$ and $\nu=1:$
\begin{equation}\label{4.10-0}
w_{0,0}(x) =-K_{0,2}(\pi,\pi-x;q_0), \quad w_{1,0}(x) =P_0(\pi,\pi-x;q_0),
\end{equation}
\begin{equation}\label{4.10-1}
w_{0,1}(x) =K_{1,2}(\pi,\pi-x;p), \quad w_{1,1}(x) =P_1(\pi,\pi-x;p),
\end{equation}
respectively. Here and below, in order to emphasize dependence of a certain function $F(x_1,x_2)$ on some
function $f(x),$ sometimes we write $F(x_1,x_2;f).$

According to (\ref{4.1})--(\ref{4.3}) and (\ref{4.10.1}), the subsystem~(\ref{4.10-0}) is nonlinear when
$a_0\in[\pi/3,\pi/2),$ while the subsystem (\ref{4.10-1}) is always linear because $a_1\ge\pi/2.$

Consider first the linear subsystem (\ref{4.10-1}). By virtue of (\ref{4.1}) and (\ref{4.10.1}), we get
$$
K_{1,1}(x,t)= \frac12 p\Big(x-\frac{t-a_1}2\Big), \quad K_{1,2}(x,t)= \frac14\Big(p\Big(\frac{t+a_1}2\Big)
-p\Big(x-\frac{t-a_1}2\Big)\Big),
$$
$$
P_1(x,t)= \frac14\Big(p\Big(\frac{t+a_1}2\Big) +p\Big(x-\frac{t-a_1}2\Big)\Big), \quad a_1<t< x<\pi.
$$
Thus, the subsystem (\ref{4.10-1}) is equivalent to the system
$$
w_{j,1}(x)= \frac14\Big(p\Big(\frac{\pi+a_1-x}2\Big) -(-1)^j p\Big(\frac{\pi+a_1+x}2\Big)\Big), \quad 0< x<
\pi-a_1, \quad j=0,1.
$$
Solving this linear system, we get
$$
p\Big(\frac{\pi+a_1-x}2\Big)=2(w_{0,1}+w_{1,1})(x), \quad p\Big(\frac{\pi+a_1+x}2\Big)=2(w_{1,1}-w_{0,1})(x),
\quad x\in(0,\pi-a_1),
$$
or, after changing the variable, we have
\begin{equation}\label{4.17}
p(x)=2\left\{\begin{array}{l}\displaystyle (w_{1,1}+w_{0,1})(\pi+a_1-2x), \quad
a_1<x<\frac{a_1+\pi}2,\\[3mm]
\displaystyle (w_{1,1}-w_{0,1})(2x-\pi-a_1), \quad \frac{a_1+\pi}2<x<\pi.
\end{array}\right.
\end{equation}

Thus, we arrive at the following theorem.

\medskip
{\bf Theorem 5. }{\it Let $a_1\in[\pi/2,\pi).$ Then for any functions $w_{0,1}(x),w_{1,1}(x)\in
L_2(0,\pi-a_1)$ the linear subsystem (\ref{4.10-1}) has a unique solution $p(x)\in L_2(a_1,\pi),$ which can
be constructed by formula (\ref{4.17})}. Moreover,
\begin{equation}\label{4.17.1}
\int\limits_{a_1}^\pi p(x)\,dx=2\int\limits_0^{\pi-a_1} w_{1,1}(x)\,dx, \quad \int\limits_{a_1}^\pi
xp(x)\,dx=\frac{\pi+a_1}2\int\limits_{a_1}^\pi p(x)\,dx -\int\limits_0^{\pi-a_1} xw_{0,1}(x)\,dx.
\end{equation}

\medskip
{\it Proof.} It remains to prove (\ref{4.17.1}). Changing the integration variable, one can easily obtain the
following relations for any integrable function $f(x):$
\begin{equation}\label{4.17.2}
\int\limits_a^\frac{\pi+a}2 f(\pi+a-2x)\,dx=\frac12\int\limits_0^{\pi-a} f(x)\,dx, \quad
\int\limits_\frac{\pi+a}2^\pi f(2x-\pi-a)\,dx=\frac12\int\limits_0^{\pi-a} f(x)\,dx,
\end{equation}
which along with (\ref{4.17}) give the first relation in (\ref{4.17.1}). Analogously, using (\ref{4.17}) and
the relations
$$
\int\limits_{a}^\frac{\pi+a}2 xf(\pi+a-2x)\,dx=\int\limits_0^{\pi-a} \frac{\pi+a-x}4f(x)\,dx, \;
\int\limits_\frac{\pi+a}2^\pi xf(2x-\pi-a)\,dx=\int\limits_0^{\pi-a} \frac{\pi+a+x}4f(x)\,dx,
$$
one can obtain the second identity in (\ref{4.17.1}). $\hfill\Box$

\medskip
Further, differentiating (\ref{4.2}) and (\ref{4.3}), and taking (\ref{4.10.1}) into account, we get
\begin{equation}\label{4.18}
K_{0,l}(x,t)= \frac12\left\{\begin{array}{r}\displaystyle q_0\Big(x-\frac{t-a_0}2\Big)
+\frac{\partial}{\partial x}A(x,t), \quad l=1,\\[3mm]
\displaystyle-\frac12\Big(q_0\Big(\frac{t+a_0}2\Big) +q_0\Big(x-\frac{t-a_0}2\Big)\Big)
+\frac{\partial}{\partial t}A(x,t), \quad l=2,
\end{array}\right.
\end{equation}
where
\begin{equation}\label{4.19}
\frac{\partial}{\partial x}A(x,t)
=\left\{\begin{array}{l} 0, \qquad a_0\le t \le\min\{x,2a_0\}, \;\; x\le\pi,\\[3mm]
\displaystyle 2\int\limits_{x-\frac{t}2+a_0}^xq_0(\tau)K_0(\tau-a_0,2(\tau-x)+t-a_0)\,d\tau, \quad 2a_0< t\le
x\le\pi,
\end{array}\right.
\end{equation}
and
\begin{equation}\label{4.20}
\frac{\partial}{\partial t}A(x,t)=
\left\{\begin{array}{l} 0, \qquad a_0\le t \le\min\{x,2a_0\}, \;\; x\le\pi,\\[3mm]
\displaystyle \int\limits_t^xq_0(\tau) K_0(\tau-a_0,t-a_0)\,d\tau\\[3mm]
\displaystyle\qquad-\int\limits_{x-\frac{t}2+a_0}^xq_0(\tau)K_0(\tau-a_0,2(\tau-x)+t-a_0)\,d\tau\\[3mm]
\qquad\qquad\quad \displaystyle -\int\limits_{\frac{t}2+a_0}^tq_0(\tau) K_0(\tau-a_0,2\tau-t-a_0)\,d\tau,
\quad 2a_0< t\le x\le\pi.
\end{array}\right.
\end{equation}
By virtue of (\ref{4.10.1}) and (\ref{4.18}), we get
\begin{equation}\label{4.21}
P_0(x,t)= \frac14\Big(q_0\Big(x-\frac{t-a_0}2\Big) -q_0\Big(\frac{t+a_0}2\Big)\Big) +B(x,t),
\end{equation}
where
$$
B(x,t)=\frac12\Big(\frac{\partial}{\partial x}A(x,t) +\frac{\partial}{\partial t}A(x,t)\Big).
$$
Then, summing up (\ref{4.19}) and (\ref{4.20}) and dividing by $2,$ we arrive at
\begin{equation}\label{4.22}
B(x,t)=
\frac12\left\{\begin{array}{l} 0, \qquad a_0\le t \le\min\{x,2a_0\}, \;\; x\le\pi,\\[3mm]
\displaystyle \int\limits_t^xq_0(\tau) K_0(\tau-a_0,t-a_0)\,d\tau\\[3mm]
\displaystyle\qquad+\int\limits_{x-\frac{t}2+a_0}^xq_0(\tau)K_0(\tau-a_0,2(\tau-x)+t-a_0)\,d\tau\\[3mm]
\qquad\qquad\quad \displaystyle -\int\limits_{\frac{t}2+a_0}^tq_0(\tau) K_0(\tau-a_0,2\tau-t-a_0)\,d\tau,
\quad 2a_0< t\le x\le\pi.
\end{array}\right.
\end{equation}
According to formula (\ref{4.18}) for $l=2$ along with (\ref{4.21}), the subsystem (\ref{4.10-0}) takes the
form
\begin{equation}\label{4.23}
w_{j,0}(x)= \frac14\Big(q_0\Big(\frac{\pi+x+a_0}2\Big) +(-1)^j q_0\Big(\frac{\pi-x+a_0}2\Big)\Big) +u_j(x),
\;\; 0\le x\le \pi-a_0, \;\; j=0,1,
\end{equation}
where
\begin{equation}\label{4.23.1}
u_0(x)=-\frac12\frac\partial{\partial t}A(\pi,t)\Big|_{t=\pi-x} =\frac12\frac{d}{dx}A(\pi,\pi-x), \quad
u_1(x)=B(\pi,\pi-x).
\end{equation}
Thus, by virtue of (\ref{4.20}) and (\ref{4.22}), we get
\begin{equation}\label{4.24}
u_j(x)= \frac12\left\{\begin{array}{l}
\displaystyle (-1)^{j+1}\int\limits_{\pi-x}^\pi q_0(\tau) K_0(\tau-a_0,\pi-x-a_0)\,d\tau\\[3mm]
\displaystyle\qquad\qquad\qquad
+\int\limits_{\frac{\pi+x}2+a_0}^\pi q_0(\tau)K_0(\tau-a_0,2\tau-x-\pi-a_0)\,d\tau\\[3mm]
\quad \displaystyle +(-1)^j\int\limits_{\frac{\pi-x}2+a_0}^{\pi-x}q_0(\tau)
K_0(\tau-a_0,2\tau+x-\pi-a_0)\,d\tau,
\quad 0\le x<\pi-2a_0,\\[3mm]
0, \qquad \pi-2a_0\le x\le \pi-a_0,
\end{array}\right.
\end{equation}
Transforming (\ref{4.23}) and taking (\ref{4.24}) into account, we obtain
$$
2(w_{0,0}+(-1)^jw_{1,0})(x)=q_0\Big(\frac{\pi+a_0+(-1)^jx}2\Big)+\left\{\begin{array}{l}2(u_0+(-1)^ju_1)(x),
\quad 0\le x<\pi-2a_0,\\[3mm]
0, \quad \pi-2a_0\le x\le \pi-a_0,
\end{array}\right.
$$
for $j=0,1.$ After changing the variables, we get
\begin{equation}\label{4.25}
2(w_{0,0}-w_{1,0})(\pi+a_0-2x)=\left\{\begin{array}{l}\displaystyle q_0(x),
\quad a_0\le x\le \frac{3a_0}2,\\[3mm]
\displaystyle q_0(x)+2v(x), \quad \frac{3a_0}2<x\le \frac{a_0+\pi}2,
\end{array}\right.
\end{equation}
\begin{equation}\label{4.26}
2(w_{0,0}+w_{1,0})(2x-\pi-a_0)=\left\{\begin{array}{l}\displaystyle q_0(x)+2v(x),
\quad \frac{a_0+\pi}2\le x<\pi-\frac{a_0}2,\\[3mm]
\displaystyle q_0(x), \quad \pi-\frac{a_0}2\le x\le \pi,
\end{array}\right.
\end{equation}
where (note that $u_1(0)=0)$
\begin{equation}\label{4.27}
v(x)=\left\{\begin{array}{l}\displaystyle (u_0-u_1)(\pi+a_0-2x),
\quad \frac{3a_0}2<x\le \frac{a_0+\pi}2,\\[3mm]
\displaystyle (u_0+u_1)(2x-\pi-a_0), \quad \frac{a_0+\pi}2<x<\pi-\frac{a_0}2.
\end{array}\right.
\end{equation}

Obviously, formulae (\ref{4.25}) and (\ref{4.26}) immediately give the solution $q_0(x)$ of the
subsystem~(\ref{4.10-0}) on $I_2:=[a_0,3a_0/2]\cup[\pi-a_0/2,\pi].$ For $a_0\ge\pi/2,$ we have
$[a_0,\pi]\subset I_2$ and, hence, the function $q_0(x)$ is completely obtained. For $a_0<\pi/2,$ the
subsystem (\ref{4.10-0}) becomes nonlinear, and its solvability on $I_3:=(3a_0/2,\pi-a_0/2)$ is conditioned
by the following lemma.

\medskip
{\bf Lemma 4. }{\it Let $v(x)$ be determined on the interval $I_3$ by formula (\ref{4.27}) with $u_0(x)$ and
$u_1(x)$ constructed in (\ref{4.24}). Then $v|_D$ does not depend on $q_0|_{I_3}$ if and only if
$D\subset[\pi-a_0,2a_0].$

Here and below $f|_{\cal S}$ denotes the restriction of the function $f$ to the set ${\cal S}.$}

\medskip
{\it Proof.} Substituting (\ref{4.24}) into (\ref{4.27}), we obtain the formulae
$$
v(x)= \int\limits_{x+\frac{a_0}2}^{2x-a_0} q_0(\tau) K_0(\tau-a_0,2(\tau-x))\,d\tau
\qquad\qquad\qquad\qquad\qquad\qquad
$$
\begin{equation}\label{4.28.1}
\qquad\qquad\qquad -\int\limits_{2x-a_0}^\pi q_0(\tau)K_0(\tau-a_0,2(x-a_0))\,d\tau, \quad
\frac{3a_0}2<x\le\frac{a_0+\pi}2,
\end{equation}
and
\begin{equation}\label{4.28.2}
v(x)= \int\limits_{x+\frac{a_0}2}^\pi q_0(\tau) K_0(\tau-a_0,2(\tau-x))\,d\tau, \quad
\frac{a_0+\pi}2<x<\pi-\frac{a_0}2.
\end{equation}
where, according to Remark~1, the function $K_0(\tau-a_0,\,\cdot\,)$ is determined by~(\ref{4.3.1}), i.e.
$$
K_0(\tau-a_0,2(\tau-x))=-K_0(\tau-a_0,2(x-a_0))
=\frac12\int\limits_{\tau-x+\frac{a_0}2}^{x-\frac{a_0}2}q_0(\zeta)\,d\zeta.
$$
Substituting this into (\ref{4.28.1}) and (\ref{4.28.2}), we arrive at
$$
v(x)=\frac12\int\limits_{x+\frac{a_0}2}^{2x-a_0} q_0(\tau)\,d\tau
\int\limits_{\tau-x+\frac{a_0}2}^{x-\frac{a_0}2}q_0(\zeta)\,d\zeta - \frac12\int\limits_{2x-a_0}^\pi
q_0(\tau)\,d\tau \int\limits_{x-\frac{a_0}2}^{\tau-x+\frac{a_0}2}q_0(\zeta)\,d\zeta  \qquad\qquad\qquad\qquad
$$
\begin{equation}\label{4.30.1}
\qquad\qquad\qquad\qquad =\frac12\int\limits_{x+\frac{a_0}2}^\pi q_0(\tau)\,d\tau
\int\limits_{\tau-x+\frac{a_0}2}^{x-\frac{a_0}2}q_0(\zeta)\,d\zeta, \quad \frac{3a_0}2<x<\pi-\frac{a_0}2.
\end{equation}
Thus, the function $v(x)$ does not depend on $q_0|_{[\pi-a_0,2a_0]}.$ Moreover, according to the first
representation in (\ref{4.30.1}), it depends on $q_0|_{(3a_0/2,\pi-a_0)}$ if and only if
$$
x-\frac{a_0}2>\frac{3a_0}2 \quad {\rm or} \quad \pi-x+\frac{a_0}2>\frac{3a_0}2,
$$
i.e. $x\in I_3\setminus[\pi-a_0,2a_0].$ Analogously, $v(x)$ depends on $q_0|_{(2a_0,\pi-a_0/2)}$ if and only
if $x<\pi-a_0.$ Hence, $v(x)$ is independent of $q_0|_{I_3}$ if and only if $x\in[\pi-a_0,2a_0].$
$\hfill\Box$

\medskip
Lemma~4 along with formulae (\ref{4.25}) and (\ref{4.26}) guaranties solvability of the subsystem
(\ref{4.10-0}) on the set $I_1=(a_0,3a_0/2)\cup (\pi-a_0,2a_0)\cup (\pi-a_0/2,\pi).$ The following corollary
gives a condition of the solvability on the entire interval $(a_0,\pi)$

\medskip
{\bf Corollary 1. }{\it $v|_{I_3}$ does not depend on $q_0|_{I_3}$ if and only if $a_0\ge2\pi/5.$}

\medskip
{\it Proof.} According to Lemma~4, $v|_{I_3}$ does not depend on $q_0|_{I_3}$ if and only if
$I_3\subset[\pi-a_0,2a_0],$ which, in turn, is equivalent to $a_0\ge2\pi/5.$ $\hfill\Box$

\medskip
According to Corollary~1, formulae (\ref{4.25}), (\ref{4.26}) and (\ref{4.30.1}) give the representation
$$
v(x)=2\int\limits_{x+\frac{a_0}2}^\pi (w_{0,0}+w_{1,0})(2\tau-\pi-a_0)\,d\tau
\int\limits_{\tau-x+\frac{a_0}2}^{x-\frac{a_0}2}(w_{0,0}-w_{1,0})(\pi+a_0-2\zeta)\,d\zeta \qquad\qquad\qquad
$$
\begin{equation}\label{4.30.2}
\qquad\qquad\qquad =\frac12\int\limits_{2x-\pi}^{\pi-a_0} (w_{0,0}+w_{1,0})(\tau)\,d\tau
\int\limits_{\pi+2a_0-2x}^{2x-a_0-\tau}(w_{0,0}-w_{1,0})(\zeta)\,d\zeta, \quad
\frac{3a_0}2<x<\frac{a_0+\pi}2,
\end{equation}
as soon as $a_0\in[2\pi/5,\pi/2).$ Thus, we arrive at the following theorem.

\medskip
{\bf Theorem 6. }{\it Let $a_0\in[2\pi/5,\pi).$ Then for any functions $w_{0,0}(x),w_{1,0}(x)\in
L_2(0,\pi-a_0)$ subsystem (\ref{4.10-0}) has a unique solution $q_0(x)\in L_2(a_0,\pi),$ which can be
constructed by the formula
\begin{equation}\label{4.31}
q_0(x)=2\left\{\begin{array}{lc} (w_{0,0}-w_{1,0})(\pi+a_0-2x), &\displaystyle a_0<x<\frac{3a_0}2,\\[3mm]
(w_{0,0}-w_{1,0})(\pi+a_0-2x)-v(x), &\displaystyle \frac{3a_0}2<x<\frac{a_0+\pi}2,\\[3mm]
(w_{0,0}+w_{1,0})(2x-\pi-a_0)-v(x), &\displaystyle \frac{a_0+\pi}2<x<\pi-\frac{a_0}2,\\[3mm]
(w_{0,0}+w_{1,0})(2x-\pi-a_0), &\displaystyle \pi-\frac{a_0}2<x<\pi,
\end{array}\right.
\end{equation}
where the function $v(x)$ is determined by formula (\ref{4.30.2}).  Moreover,
\begin{equation}\label{4.31.1}
\int\limits_{a_0}^\pi q_0(x)\,dx=2\int\limits_0^{\pi-a_0} w_{0,0}(x)\,dx.
\end{equation}}

{\it Proof.} It remains to prove (\ref{4.31.1}). Indeed, integrating (\ref{4.31}) and using (\ref{4.17.2}),
we get
$$
\int\limits_{a_0}^\pi q_0(x)\,dx=2\int\limits_0^{\pi-a_0} w_{0,0}(x)\,dx
-2\int\limits_\frac{3a_0}2^{\pi-\frac{a_0}2} v(x)\,dx,
$$
where, according to (\ref{4.27}), we have
$$
\int\limits_\frac{3a_0}2^{\pi-\frac{a_0}2} v(x)\,dx=\int\limits_\frac{3a_0}2^\frac{\pi+a_0}2
(u_0-u_1)(\pi+a_0-2x)\,dx +\int\limits_\frac{\pi+a_0}2^{\pi-\frac{a_0}2}
(u_0+u_1)(2x-\pi-a_0)\,dx=\int\limits_0^{\pi-2a_0} u_0(x)\,dx.
$$
Finally, using (\ref{4.23.1}) and the last two equalities in (\ref{4.4}), we calculate
$$
\int\limits_0^{\pi-2a_0} u_0(x)\,dx=\frac{A(\pi,2a_0)-A(\pi,\pi)}2=0,
$$
which finishes the proof. $\hfill\Box$
\\

{\large\bf 4. Proof of the uniqueness theorems}
\\

Let us begin with the following assertion.

\medskip
{\bf Lemma 5. }{\it Let $a_0,a_1\in(0,\pi).$ Then specification of any pair of sequences $\{\rho_{n,j}\}_{
n\in{\mathbb Z}_j},$ $j=0,1,$ of the form~(\ref{2.1}) uniquely determines the values $\alpha_0,$ $\alpha_1$
and $\omega.$}

\medskip
{\it Proof.} According to (\ref{2.1}), we arrive at the asymptotic formulae
\begin{equation}\label{4.34.2}
\omega\cos na_0=\frac{\gamma_{n,0}+\gamma_{-n,0}}2 +o(1), \quad
\alpha_j\sin\rho_{n,j}^0a_1=\frac{\gamma_{n,j}-\gamma_{j-n,j}}2 +o(1), \quad |n|\to\infty,
\end{equation}
where we denoted
\begin{equation}\label{4.34.1}
\gamma_{n,j}:=\pi n(\rho_{n,j}-\rho_{n,j}^0),  \quad n\in{\mathbb Z}_j, \quad j=0,1.
\end{equation}
By virtue of Lemma~3.3 in \cite{SS19}, the sequences $\{\cos na_0\}_{n\ge1}$ and $\{\sin
\rho_{n,j}^0a_1\}_{n\ge1},$ $j=0,1,$ do not converge, i.e. each of them has at least two different partial
limits. Choose increasing sequences of natural numbers $\{m_{k,l}\},$ $l=\overline{1,3},$ so that
\begin{equation}\label{4.34.4}
r_1:=\lim_{k\to\infty}\cos m_{k,1}a_0\ne0, \quad r_l:=\lim_{k\to\infty}\sin\rho_{m_{k,l},l-2}^0a_1\ne0, \quad
l=2,3.
\end{equation}
According to (\ref{4.34.2}) and (\ref{4.34.4}), we calculate the values $\omega,$ $\alpha_0$ and $\alpha_1$
by the formulae
\begin{equation}\label{4.34.5}
\omega=\lim_{k\to\infty}\frac{\gamma_{m_{k,1},0}+\gamma_{-m_{k,1},0}}{2r_1}, \quad
\alpha_j=\lim_{k\to\infty}\frac{\gamma_{m_{k,j+2},j}-\gamma_{j-m_{k,j+2},j}}{2r_{j+2}}, \quad j=0,1,
\end{equation}
which finishes the proof. $\hfill\Box$

\medskip
{\it Proof of Theorem~2.} According to Lemmas~3 and~5, under the hypothesis of the theorem, the functions
$\Delta_j(\rho),\,j=0,1,$ as well as the numbers $\alpha_0,$ $\alpha_1,$ and $\omega$ are determine uniquely.
Then, by virtue of Lemma~2, the functions $w_{j,\nu}(x),j,\nu=0,1,$ are uniquely determined too. Thus, by
Theorem~5 and Lemma~4, the function $p(x)=q_1'(x)$ is uniquely determined a.e. on the interval $(a_1,\pi),$
while $q_0(x)$ is so on $I_1.$ Finally, taking~(\ref{2.2}) into account, we get
\begin{equation}\label{4.39.1}
q_1(x)=\alpha_0+\alpha_1 +\int\limits_{a_1}^x p(t)\,dt, \quad a_1\le x\le\pi,
\end{equation}
which finishes the proof. $\hfill\Box$

\medskip
{\it Proof of Theorem~3.} Following the proof of Theorem~2 in \cite{DV19}, we denote
$$
p_1:=q_0|_{(a_0,\pi-a_0)}, \quad p_2:=q_0|_{(2a_0,\pi)}, \quad
R_1(x,t):=\int\limits_{t-x+\frac{a_0}2}^{x-\frac{a_0}2} p_1(\tau)\,d\tau,\quad R_2(x):=\int\limits_x^\pi
p_2(t)\,dt.
$$
Then, by virtue of (\ref{4.30.1}) and (\ref{4.31}) we have the relation
\begin{equation}\label{4.39.2}
F(x) =q_0(x)+\int\limits_{x+\frac{a_0}2}^\pi p_2(t)\,dt
\int\limits_{t-x+\frac{a_0}2}^{x-\frac{a_0}2}p_1(\tau)\,d\tau, \quad \frac{3a_0}2<x<\pi-\frac{a_0}2,
\end{equation}
where
$$
F(x)=2\left\{\begin{array}{lc}
(w_{0,0}-w_{1,0})(\pi+a_0-2x), &\displaystyle \frac{3a_0}2<x<\frac{a_0+\pi}2,\\[3mm]
(w_{0,0}+w_{1,0})(2x-\pi-a_0), &\displaystyle \frac{a_0+\pi}2<x<\pi-\frac{a_0}2.
\end{array}\right.
$$
After changing the order of integration, relation (\ref{4.39.2}) on the target intervals takes the forms:
\begin{equation}\label{4.39.3}
F_1(x) =p_1(x)+\int\limits_{x+\frac{a_0}2}^{\pi-\frac{a_0}2}R_1(x,t) p_2(t)\,dt -
\int\limits_{\frac{3a_0}2}^{\pi-x+\frac{a_0}2}R_2\Big(x+t-\frac{a_0}2\Big)p_1(t)\,dt, \;\;
\frac{3a_0}2<x<\pi-a_0,
\end{equation}
where
$$
F_1(x) :=F(x)+ \int\limits_{\pi-x}^{\frac{3a_0}2}R_2\Big(x+t-\frac{a_0}2\Big)p_1(t)\,dt
+R_2\Big(\pi-\frac{a_0}2\Big)\int\limits_{x-\frac{a_0}2}^{\pi-x}p_1(t) \,dt, \;\; \frac{3a_0}2<x<\pi-a_0,
$$
and
\begin{equation}\label{4.39.4}
F_2(x) =p_2(x)+R_2\Big(x+\frac{a_2}2\Big)\int\limits_\frac{3a_0}2^{x-\frac{a_0}2}p_1(\tau) \,d\tau, \quad
2a_0<x<\pi-\frac{a_0}2,
\end{equation}
where
$$
F_2(x) :=F(x)- \int\limits_{x+\frac{a_0}2}^\pi p_2(t) \,dt \int\limits_{t-x+\frac{a_0}2}^\frac{3a_0}2
p_1(\tau)\,d\tau, \quad 2a_0<x<\pi-\frac{a_0}2.
$$
Note that the functions $F_1(x),$ $F_2(x),$ $R_1(x,t)$ and $R_2(x)$ involved in (\ref{4.39.3}) and
(\ref{4.39.4}), being dependent only on $q|_{I_2},$ are already known. Substituting (\ref{4.39.4}) into
(\ref{4.39.3}) and changing the order of integration, we obtain the integral equation
$$
F_3(x) =p_1(x) -R_3(x,x)\int\limits_\frac{3a_0}2^x p_1(t)\,dt -\int\limits_x^{\pi-a_0} R_3(x,t)p_1(t)\,dt
\qquad\qquad\qquad\;
$$
\begin{equation}\label{4.39.5}
\qquad\qquad\quad -\int\limits_{\frac{3a_0}2}^{\pi-x+\frac{a_0}2} R_2\Big(x+t-\frac{a_0}2\Big)p_1(t)\,dt,
\quad \frac{3a_0}2<x<\pi-a_0,
\end{equation}
where
$$
F_3(x) =F_1(x)- \int\limits_{x+\frac{a_0}2}^{\pi-\frac{a_0}2}R_1(x,t)F_2(t)\,dt, \quad R_3(x,t)=
\int\limits_{t+\frac{a_0}2}^{\pi-\frac{a_0}2}R_1(x,\tau)R_2\Big(\tau-\frac{a_0}2\Big)\,d\tau.
$$
Since, according to the hypothesis of the theorem, the value
$$
\omega_1:=\omega_0+\int\limits_\frac{3a_0}2^{\frac\pi2+\frac{a_0}4}q_0(x)\,dx
=\int\limits_\frac{3a_0}2^{\pi-a_0}p_1(x)\,dx
$$
is known, equation (\ref{4.39.5}) takes the form
\begin{equation}\label{4.39.6}
F_4(x) =p_1(x) +\int\limits_x^{\pi-a}R_4(x,t) p_1(t)\,dt -\int\limits_{\frac{3a_0}2}^{\pi-x+\frac{a_0}2}
R_2\Big(x+t-\frac{a_0}2\Big)p_1(t)\,dt, \quad \frac{3a_0}2<x<\pi-a_0,
\end{equation}
where both functions $F_4(x)=F_3(x)+\omega_1R_3(x,x)$ and $R_4(x,t)=R_3(x,x)-R_3(x,t)$ are known too. Taking
into account that the function $p_1(x)$ on the interval $(3a_0/2,\pi/2+a_0/4)$ is given by the hypothesis, we
find it on $(\pi/2+a_0/4,\pi-a_0)$ by solving there equation (\ref{4.39.6}), in which the second integral
becomes known. Finally, substituting the completely found function $p_1(x)$ into relation (\ref{4.39.4}), we
find $p_2(x)$ on $(2a_0,\pi-a_0/2),$ which finishes the proof. $\hfill\Box$
\\

{\large\bf 5. Other representations for the infinite products}
\\

The results of this section are valid for any $a_0,a_1\in[0,\pi].$

In Section~3, we proved, in particular, that any functions $\Delta_0(\lambda)$ and $\Delta_1(\lambda)$ of the
forms described in (\ref{4.11})--(\ref{4.12.1.2}) have infinitely many zeros obeying (\ref{2.1}) for the
corresponding $j\in\{0,1\}.$ Moreover, these functions are determined by their zeros uniquely by formula
(\ref{4.13}). However, the functions $\Delta_0(\lambda)$ and $\Delta_1(\lambda)$ constructed by (\ref{4.13})
with arbitrary sequences of complex numbers of the form~(\ref{2.1}), generally speaking, do not have the
forms as in (\ref{4.11})--(\ref{4.12.1.2}). This fact is connected, in particular, with excessiveness of the
input data of Inverse Problem~1 for recovering the functions~$w_{j,\nu}(x)$ in (\ref{4.11}) and~(\ref{4.12}).
Nevertheless, the following lemma holds.

\medskip
{\bf Lemma 6. }{\it Let $j\in\{0,1\}.$ Then for any sequence of complex numbers $\{\rho_{n,j}\}_{n\in{\mathbb
Z}_j}$ obeying~(\ref{2.1}), the function $\Delta_j(\rho)$ constructed by formula (\ref{4.13}) has the form
\begin{equation}\label{4.11.1}
\Delta_0(\rho)=\frac{\sin\rho\pi}\rho-\omega\frac{\cos\rho(\pi-a_0)}{\rho^2}
+\alpha_0\frac{\sin\rho(\pi-a_1)}{\rho^2} +\gamma_0\frac{\sin\rho\pi}{\rho^2}
+\sum_{\nu=0}^1\int\limits_0^\pi w_{0,\nu}(x) \frac{c_\nu(\rho x)}{\rho^2}\,dx
\end{equation}
for $j=0,$ and
\begin{equation}\label{4.12.1}
\Delta_1(\rho)=\cos\rho\pi+\omega\frac{\sin\rho(\pi-a_0)}\rho +\alpha_1\frac{\cos\rho(\pi-a_1)}\rho
+\gamma_1\frac{\cos\rho\pi}\rho +\sum_{\nu=0}^1\int\limits_0^\pi w_{1,\nu}(x) \frac{c_{1-\nu}(\rho
x)}\rho\,dx
\end{equation}
for $j=1.$ Here $w_{j,\nu}(x)\in L_2(0,\pi)$ for $j,\nu=0,1,$ and
\begin{equation}\label{4.12.1.1}
\int\limits_0^\pi w_{0,0}(x)\,dx=\omega, \;\; \int\limits_0^\pi
xw_{0,1}(x)\,dx=\alpha_0(a_1-\pi)-\gamma_0\pi,\;\; \int\limits_0^\pi w_{1,1}(x)\,dx=-\alpha_1 -\gamma_1.
\end{equation}}

Before proceeding directly to the proof of Lemma~6, we establish the following auxiliary assertion giving
some important subtle estimates that will be required for the proof.

\medskip
{\bf Proposition 1. }{\it Put $h_{n,j}:=\beta_0\cos\rho_{n,j}^0a_0 +\beta_1\sin\rho_{n,j}^0a_1+\varkappa_n,$
$|n|\in{\mathbb N},$ $\beta_j\in{\mathbb C},$ $j=0,1.$ Then
$$
a_{n,j}:=\sum_{k\ne0,n}\frac{h_{k,j}}{k(n-k)} =O\Big(\frac1n\Big), \quad |n|\to\infty.
$$}

{\it Proof.} We have
$$
a_{n,j}=\frac1n\lim_{N\to\infty}\sum_{{k\ne0,n}\atop{k=-N}}^Nh_{k,j}\Big(\frac1k +\frac1{n-k}\Big)
=\frac1n\lim_{N\to\infty}\Big(\sum_{{k\ne0,n}\atop{k=-N}}^N\frac{h_{k,j}}k
-\sum_{{k\ne-n,0}\atop{k=-N-n}}^{N-n}\frac{h_{n+k,j}}k\Big).
$$
One can easily calculate
$$
a_{n,j}=\frac1n \sum_{{k\ne|n|}\atop{k=1}}^\infty\frac{h_{k,n,j}}k-\frac{h_{-n,j}+h_{2n,j}}{n^2}, \quad
h_{k,n,j}=h_{k,j}-h_{n+k,j}+h_{n-k,j}-h_{-k,j}.
$$
Denote $\varkappa_{k,n}:= \varkappa_k-\varkappa_{n+k}+\varkappa_{n-k}-\varkappa_{-k}.$ It remains to note
that, since
$$
h_{k,n,j} = 4\beta_0\sin\frac{na_0}2\cos\frac{(n-j)a_0}2\sin ka_0
+4\beta_1\sin\frac{na_1}2\sin\frac{(n-j)a_1}2\sin ka_1 +\varkappa_{k,n}
$$
and the series $\sum_{k=1}^\infty\frac{\sin ka}k$ converges for any $a,$ the series
$\sum_{k=1}^\infty\frac{h_{k,n,j}}k,$ $j=0,1,$ are convergent too, and their sums are uniformly bounded with
respect to $n.$ $\hfill\Box$

\medskip
{\it Proof of Lemma~6.} Fix $j\in\{0,1\}.$ Let us show first that $\{\theta_j(\rho_{n,j}^0)\}_{n\in{\mathbb
Z}_j}\in l_2,$ where $\theta_0(\rho)$ and $\theta_1(\rho)$ are determined by (\ref{4.33}) and (\ref{4.34}),
respectively. It is easy to see that
\begin{equation}\label{4.12.1.1.1}
\theta_j(\rho_{n,j}^0) =(\rho_{n,j}^0)^{2-j}\Delta_j(\rho_{n,j}^0) +(-1)^n\omega_{n,j}, \quad
\omega_{n,j}=\omega\cos \rho_{n,j}^0a_0 +\alpha_j\sin \rho_{n,j}^0a_1.
\end{equation}
By virtue of (\ref{4.13}) and (\ref{4.15}), we obtain
$$
\rho^{2-j}\Delta_j(\rho)=\rho(\rho_{n,j}-\rho)\frac{c_{1-j}(\rho\pi)}{\rho_{n,j}^0-\rho} \prod_{k\in{\mathbb
Z}_j\setminus\{n\}} \frac{\rho_{k,j}-\rho}{\rho_{k,j}^0-\rho}.
$$
Thus, having put
\begin{equation}\label{4.12.1.1.2}
\varepsilon_{n,j}:=\rho_{n,j}-\rho_{n,j}^0=\frac{\omega_{n,j}}{\pi n}+\frac{\varkappa_n}n,
\end{equation}
we get
$$
(\rho_{n,j}^0)^{2-j}\Delta_j(\rho_{n,j}^0)=(-1)^{n+1}\pi\rho_{n,j}^0\varepsilon_{n,j}b_{n,j},  \quad
b_{n,j}=\prod_{k\in{\mathbb Z}_j\setminus\{n\}}\Big(1+\frac{\varepsilon_{k,j}}{k-n}\Big), \quad n\in{\mathbb
Z}_j.
$$
Since $\pi\rho_{n,j}^0\varepsilon_{n,j}=\omega_{n,j} +\varkappa_n,$ we get $\theta_j(\rho_{n,j}^0)
=(-1)^n(1-b_{n,j})\omega_{n,j}+\varkappa_n.$ Thus, we need to prove that $\{1-b_{n,j}\}\in l_2.$ For this
purpose, we choose $N\in{\mathbb N}$ so that $|\varepsilon_{n,j}|\le1/2$ as soon as $|n|\ge N,$ and represent
$b_{n,j}$ in the form $b_{n,j}=b_{n,j}^{(1)}b_{n,j}^{(2)},$ where
$$
b_{n,j}^{(1)}=\prod_{{k\ne n}\atop{1-j\le|k|<N}}\Big(1+\frac{\varepsilon_{k,j}}{k-n}\Big)
=1+O\Big(\frac1n\Big), \;\; |n|\to\infty, \quad b_{n,j}^{(2)}=\prod_{{k\ne n}\atop{|k|\ge
N}}\Big(1+\frac{\varepsilon_{k,j}}{k-n}\Big).
$$
Our choice of $N$ allows one to represent
$$
b_{n,j}^{(2)}=\exp\Big(\sum_{{k\ne n}\atop{|k|\ge N}}\ln\Big(1+\frac{\varepsilon_{k,j}}{k-n}\Big)\Big)
=\exp\Big(\sum_{{k\ne n}\atop{|k|\ge
N}}\sum_{\nu=0}^\infty\frac{(-1)^\nu}{\nu+1}\Big(\frac{\varepsilon_{k,j}}{k-n}\Big)^{\nu+1}\Big).
$$
Therefore, we have
$$
|b_{n,j}^{(2)}-1| \le\sum_{\nu=1}^\infty\frac{\Big( \Omega_{n,j}^{(1)}+2\Omega_{n,j}^{(2)}\Big)^\nu}{\nu!},
\quad \Omega_{n,j}^{(1)}=\Big|\sum_{{k\ne n}\atop{|k|\ge N}}\frac{\varepsilon_{k,j}}{k-n}\Big|, \quad
\Omega_{n,j}^{(2)} =\sum_{{k\ne n}\atop{|k|\ge N}}\frac{|\varepsilon_{k,j}|^2}{(k-n)^2}.
$$
By virtue of (\ref{4.12.1.1.1}) and (\ref{4.12.1.1.2}) along with Proposition~1, we have
$\Omega_{n,j}^{(1)}=O(n^{-1}),$ $|n|\to\infty.$ Let us show that $\{\Omega_{n,j}^{(2)}\}\in l_2.$ Indeed,
using the generalized Minkovskii inequality, we get the estimate
$$
\sqrt{\sum_{|n|\in{\mathbb N}}\Big(\Omega_{n,j}^{(2)}\Big)^2}\le C\sqrt{\sum_{|n|\in{\mathbb
N}}\Big(\sum_{k\ne0,n} \frac1{k^2(k-n)^2}\Big)^2} \qquad\qquad\qquad\qquad\quad\,
$$
$$
\qquad\qquad\qquad\qquad\le C\sum_{|k|\in{\mathbb N}} \sqrt{\sum_{n\ne0,k}\Big(\frac1{k^2(k-n)^2}\Big)^2} <
C\sum_{|k|\in{\mathbb N}}\frac1{k^2} \sqrt{\sum_{|n|\in{\mathbb N}}\frac1{n^4}}<\infty.
$$
Thus, we arrive at $\{\theta_j(\rho_{n,j}^0)\}_{n\in{\mathbb Z}_j}\in l_2.$ Further, since the systems of
vector-functions
$$
\Big\{[\cos \rho_{n,j}^0x,\sin \rho_{n,j}^0x]\Big\}_{n\in{\mathbb Z}},
$$
where $\rho_{0,0}^0=0,$ is an almost normalized orthogonal basis in $(L_2(0,\pi))^2,$ there exist unique
functions $w_{j,\nu}(x)\in L_2(0,\pi),$ $\nu=0,1,$ such that
$$
\theta_j(\rho_{n,j}^0)=\int\limits_0^\pi w_{j,j}(x)\cos\rho_{n,j}^0 x\,dx +\int\limits_0^\pi
w_{j,1-j}(x)\sin\rho_{n,j}^0 x\,dx, \quad n\in{\mathbb Z}.
$$
Consider the function
$$
\tilde\theta_j(\rho)=\int\limits_0^\pi w_{j,j}(x)\cos\rho x\,dx +\int\limits_0^\pi w_{j,1-j}(x)\sin\rho
x\,dx.
$$
According to (\ref{4.33}) and (\ref{4.34}), it remains to prove that
\begin{equation}\label{4.12.1.5}
\theta_j(\rho)-\tilde\theta_j(\rho)=\gamma_jc_{1-j}(\rho\pi), \quad \gamma_j\equiv const.
\end{equation}
For this purpose, we consider the entire function
\begin{equation}\label{4.12.1.3}
\Theta_j(\rho):=\frac{\theta_j(\rho)-\tilde\theta_j(\rho)}{c_{1-j}(\rho\pi)} =\Theta_{j,1}(\rho)
+\Theta_{j,2}(\rho),
\end{equation}
where
$$
\Theta_{j,1}(\rho)=\rho\Big(\frac{\rho^{1-j}\Delta_j(\rho)}{c_{1-j}(\rho\pi)}-1\Big), \quad
\Theta_{j,2}(\rho)=\frac{(-1)^j\omega c_j(\rho(\pi-a_0))
-\alpha_jc_{1-j}(\rho(\pi-a_1))-\tilde\theta_j(\rho)}{c_{1-j}(\rho\pi)}.
$$
Clearly, $\Theta_{j,2}(\rho)=O(1)$ as soon as $\rho\in G_\delta^j:=\{\rho:|\rho-\rho_{n,j}^0|\ge\delta,n\in
{\mathbb Z}\}$ for a fixed $\delta>0$ and $\rho\to\infty,$ and $\Theta_{j,2}(\rho)=o(1)$ for $|{\rm
Im}\rho|\to\infty.$ Further, dividing (\ref{4.13}) by (\ref{4.15}), we get
$$
\Theta_{j,1}(\rho)=\rho(F_{j,1}(\rho)F_{j,2}(\rho)-1),
$$
where
$$
F_{j,1}(\rho)=\prod_{1-j\le|n|<N}\Big(1+\frac{\varepsilon_{n,j}}{\rho_{n,j}^0-\rho}\Big), \quad
F_{j,2}(\rho)=\prod_{|n|\ge N}\Big(1+\frac{\varepsilon_{n,j}}{\rho_{n,j}^0-\rho}\Big),
$$
while $N$ is chosen so that $|\varepsilon_{n,j}|\le\delta/2$ as soon as $|n|\ge N.$ Hence, one can represent
$$
F_{j,2}(\rho)=\exp\Big(\sum_{|n|\ge N}\ln\Big(1+\frac{\varepsilon_{n,j}}{\rho_{n,j}^0-\rho}\Big)\Big)
=\exp\Big(\sum_{|n|\ge N}
\sum_{\nu=0}^\infty\frac{(-1)^\nu}{\nu+1}\Big(\frac{\varepsilon_{n,j}}{\rho_{n,j}^0-\rho}\Big)^{\nu+1}\Big),
\quad \rho\in G_\delta^j.
$$
Then the following estimates hold:
$$
|F_{j,2}(\rho)-1|\le \sum_{\nu=1}^\infty\frac1{\nu!}\Big(2\sum_{|n|\ge N}
\frac{|\varepsilon_{n,0}|}{|\rho_{n,j}^0-\rho|}\Big)^\nu < C \sqrt{\sum_{n\in{\mathbb Z}}
\frac1{|\rho_{n,j}^0-\rho|^2}}, \quad \rho\in G_\delta^j,
$$
where the right-hand side, as a function of $\rho,$ has the period $1.$ Thus, we arrive at the estimates
$$
F_{j,\nu}(\rho)=1+O(\rho^{\nu-2}), \;\; \rho\in G_\delta^j, \;\; \rho\to\infty, \;\;\nu=1,2, \qquad
F_{j,2}(\rho)=1+o(1), \;\; |{\rm Im}\rho|\to\infty,
$$
which imply $\Theta_{j,1}(\rho)=O(\rho)$ for $\rho\to\infty$ in $G_\delta^j,$ and $\Theta_{j,1}(\rho)
=o(\rho)$ for $|{\rm Im}\rho|\to\infty.$ Hence, we get $\gamma_j:=\Theta_j(\rho)\equiv const,$ which along
with (\ref{4.12.1.3}) implies (\ref{4.12.1.5}), i.e. (\ref{4.11.1}) and (\ref{4.12.1}) are proven.

Finally, note that relations (\ref{4.12.1.1}) follow from entireness of the functions
$\Delta_j(\rho),\,j=0,1.$ $\hfill\Box$
\\

{\large\bf 6. Solution of the inverse problem}
\\

In this section, besides our initial assumptions on $a_0$ and $a_1,$ we also assume $a_0\ge5\pi/2.$ The
preliminary work fulfilled in Sections~3 and~5 allows us to give the proof of Theorem~4.

\medskip
{\it Proof of Theorem~4.} By necessity, the asymptotics~(\ref{2.1}) is already established in Theorem~1. Let
us prove (ii). According to Lemma~3, the functions $\Delta_0(\rho)$ and $\Delta_1(\rho)$ determined by
formula (\ref{4.13}) are the characteristic functions, which, by virtue of Lemma~2, have representations
(\ref{4.11}) and (\ref{4.12}), respectively. Hence, according to (\ref{4.33}) and (\ref{4.34}), we have the
representations
$$
\theta_j(\rho)=\int\limits_0^{\pi-a_j}w_{j,j}(x)\cos\rho x\,dx
+\int\limits_0^{\pi-a_{1-j}}w_{j,1-j}(x)\sin\rho x\,dx, \quad j=0,1,
$$
which along with (\ref{4.32}) give
\begin{equation}\label{4.12.3}
g_{j,j}(\rho)=2\int\limits_0^{\pi-a_j}w_{j,j}(x)\cos\rho x\,dx, \quad
g_{j,1-j}(\rho)=2\int\limits_0^{\pi-a_{1-j}}w_{j,1-j}(x)\sin\rho x\,dx, \quad j=0,1,
\end{equation}
which, in turn, implies (ii) and finishes the proof of the necessity.

For the sufficiency, we assume that some complex sequences $\{\rho_{n,0}\}_{|n|\in{\mathbb N}}$ and
$\{\rho_{n,1}\}_{n\in{\mathbb Z}}$ obeying (i) and (ii) are given. Find the values $\alpha_0,$ $\alpha_1$ and
$\omega$ as in the proof of Lemma~5. Then, by formula (\ref{4.13}), construct the functions $\Delta_0(\rho)$
and $\Delta_1(\rho),$ which, according to Lemma~6, have representations (\ref{4.11.1}) and (\ref{4.12.1}),
respectively, with some numbers $\gamma_0$ and $\gamma_1$ and some functions $w_{\nu,j}(x)\in
L_2(0,\pi),\,j,\nu=0,1,$ obeying (\ref{4.12.1.1}). Using (\ref{4.32})--(\ref{4.34}) and (\ref{4.11.1}),
(\ref{4.12.1}), we calculate
$$
g_{0,0}(\rho)=2\int\limits_0^\pi w_{0,0}(x)\cos\rho x\,dx, \quad g_{0,1}(\rho)=2\gamma_0\sin\rho\pi
+2\int\limits_0^\pi w_{0,1}(x)\sin\rho x\,dx,
$$
$$
g_{1,0}(\rho)=2\int\limits_0^\pi w_{1,0}(x)\sin\rho x\,dx, \quad g_{1,1}(\rho)=2\gamma_1\cos\rho\pi
+2\int\limits_0^\pi w_{1,1}(x)\cos\rho x\,dx.
$$
Thus, condition (ii) implies $\gamma_0=\gamma_1=0$ and hence, by virtue of the Paley--Wiener theorem,
$w_{j,\nu}(x)=0$ a.e. on $(\pi-a_\nu,\pi)$ for $j,\nu=0,1,$ which along with (\ref{4.12.1.1}) gives
(\ref{4.12.1.2}). Therefore, the functions $\Delta_0(\lambda)$ and $\Delta_1(\lambda)$ has the forms
(\ref{4.11}) and (\ref{4.12}), respectively. By virtue of Theorems~5 and~6, the subsystems~(\ref{4.10-0}) and
(\ref{4.10-1}) with these $w_{j,\nu}(x)$ have unique solutions $q_0(x)\in L_2(a_0,\pi)$ and $p(x)\in
L_2(a_1,\pi),$ satisfying (\ref{4.31.1}) and (\ref{4.17.1}), respectively. Construct the function $q_1(x)\in
W_2^1[a_1,\pi]$ by the formula
\begin{equation}\label{4.34.7}
q_1(x)=\frac1{\pi-a_1}\int\limits_{a_1}^\pi dt\int\limits_t^\pi p(\tau)\,d\tau-\int\limits_x^\pi p(t)\,dt,
\end{equation}
which, obviously, obeys (\ref{1.1}). Thus, we constructed the boundary value problems ${\cal L}_0(q_0,q_1)$
and ${\cal L}_1(q_0,q_1).$ Let $\tilde\Delta_0(\lambda)$ and $\tilde\Delta_1(\lambda)$ be their
characteristic functions, respectively. According to Lemma~2, they have the representations
\begin{equation}\label{4.11.2}
\tilde\Delta_0(\rho)=\frac{\sin\rho\pi}\rho-\tilde\omega\frac{\cos\rho(\pi-a_0)}{\rho^2}
+\tilde\alpha_0\frac{\sin\rho(\pi-a_1)}{\rho^2} +\sum_{\nu=0}^1\int\limits_0^{\pi-a_\nu} \tilde w_{0,\nu}(x)
\frac{c_\nu(\rho x)}{\rho^2}\,dx,
\end{equation}
\begin{equation}\label{4.12.2}
\tilde\Delta_1(\rho)=\cos\rho\pi+\tilde\omega\frac{\sin\rho(\pi-a_0)}\rho
+\tilde\alpha_1\frac{\cos\rho(\pi-a_1)}\rho +\sum_{\nu=0}^1\int\limits_0^{\pi-a_\nu} \tilde w_{1,\nu}(x)
\frac{c_{1-\nu}(\rho x)}\rho\,dx,
\end{equation}
where
\begin{equation}\label{2.2-ti}
\tilde\omega=\frac12\int\limits_{a_0}^\pi q_0(x)\,dx, \quad \tilde\alpha_j=\tilde\alpha+(-1)^j\tilde\beta,
\;\; j=0,1, \quad\tilde\alpha=\frac{q_1(a_1)}2, \quad \tilde\beta=\frac{q_1(\pi)}2,
\end{equation}
\begin{equation}\label{4.10-0-ti}
\tilde w_{0,0}(x) =-K_{0,2}(\pi,\pi-x;q_0), \quad \tilde w_{1,0}(x) =P_0(\pi,\pi-x;q_0),
\end{equation}
\begin{equation}\label{4.10-1-ti}
\tilde w_{0,1}(x) =K_{1,2}(\pi,\pi-x;p), \quad \tilde w_{1,1}(x) =P_1(\pi,\pi-x;p).
\end{equation}
Comparing (\ref{4.10-0-ti}) and (\ref{4.10-1-ti}) with (\ref{4.10-0}) and (\ref{4.10-1}), respectively, we
arrive at
\begin{equation}\label{4.35.0}
\tilde w_{j,\nu}(x) = w_{j,\nu}(x), \quad j,\nu=0,1.
\end{equation}
Successively using the first equality in (\ref{2.2-ti}), identity (\ref{4.31.1}) and the first equality in
(\ref{4.12.1.2}), we get
\begin{equation}\label{4.35}
\tilde\omega=\frac12\int\limits_{a_0}^\pi q_0(x)\,dx =\int\limits_0^{\pi-a_0} w_{0,0}(x)\,dx =\omega.
\end{equation}
Further, using the first equality in (\ref{4.17.1}) along with the third one in (\ref{4.12.1.2}), we obtain
\begin{equation}\label{4.36}
\int\limits_{a_1}^\pi p(x)\,dx=2\int\limits_0^{\pi-a_1} w_{1,1}(x)\,dx=-2\alpha_1,
\end{equation}
while the second equalities in (\ref{4.12.1.2}) and in (\ref{4.17.1}) along with (\ref{4.36}) give
\begin{equation}\label{4.37}
\int\limits_{a_1}^\pi xp(x)\,dx=\frac{\pi+a_1}2\int\limits_{a_1}^\pi p(x)\,dx +(\pi-a_1)\alpha_0 =
(\pi-a_1)\alpha_0-(\pi+a_1)\alpha_1.
\end{equation}
On the other hand, successively using (\ref{2.2-ti}), (\ref{4.34.7}) and (\ref{4.36}), we get
\begin{equation}\label{4.38}
\tilde\alpha_1= \frac{q_1(a_1)-q_1(\pi)}2 =-\frac12\int\limits_{a_1}^\pi p(x)\,dx=\alpha_1,
\end{equation}
while the second and the last equalities in (\ref{2.2-ti}) along with (\ref{4.34.7}), (\ref{4.36}) and
(\ref{4.37}) imply
\begin{equation}\label{4.39}
(\pi-a_1)(\tilde\alpha_0-\tilde\alpha_1)=\int\limits_{a_1}^\pi dt\int\limits_t^\pi p(\tau)\,d\tau
=\int\limits_{a_1}^\pi (x-a_1)p(x)\,dx =(\pi-a_1)(\alpha_0-\alpha_1).
\end{equation}
By virtue of (\ref{4.38}) and (\ref{4.39}), we have $\tilde\alpha_j=\alpha_j,\,j=0,1,$ which along with
(\ref{4.11}), (\ref{4.12}),  (\ref{4.11.2}), (\ref{4.12.2}), (\ref{4.35.0}) and (\ref{4.35}) gives
$\tilde\Delta_j(\rho)\equiv\Delta_j(\rho),$ $j=0,1.$ Hence, each given sequence $\{\rho_{n,j}\}_{n\in{\mathbb
Z}_j}$ is the spectrum of the corresponding problem ${\cal L}_j(q_0,q_1),$ $j=0,1.$ $\hfill\Box$

\medskip
The Paley--Wiener theorem implies the following corollary from Theorem~4.

\medskip
{\bf Corollary 2. }{\it Arbitrary complex sequences $\{\rho_{n,0}\}_{|n|\in{\mathbb N}}$ and
$\{\rho_{n,1}\}_{n\in{\mathbb Z}}$ are the spectra of some boundary value problems ${\cal L}_0(q_0,q_1)$ and
${\cal L}_1(q_0,q_1),$ respectively, if and only if their convergence exponents are equal to~$1,$ and there
exist some $\alpha_1,\alpha_2,\,\omega\in{\mathbb C}$ such that the functions $g_{j,\nu}(\rho),$ $j,\nu=0,1,$
determined by formulae~(\ref{4.32})--(\ref{4.13}) satisfy the following conditions:
\begin{equation}\label{4.39.0}
g_{j,\nu}(x)\in L_2(-\infty,\infty), \quad |g_{j,\nu}(\rho)|\le C\exp((\pi-a_\nu)|\rho|), \quad
g_{j,\nu}(-\rho)=(-1)^{j+\nu}g_{j,\nu}(\rho).
\end{equation}}

{\it Proof.} In addition to the proof of Theorem~4, it is sufficient to note that, by virtue of the
Paley--Wiener theorem, conditions (\ref{4.39.0}) are equivalent to the representations (\ref{4.12.3}) with
some functions $w_{j,\nu}(x)\in L_2(0,\pi-a_\nu),$ $j,\nu=0,1.$ $\hfill\Box$

\medskip
The proof of Theorem~4 gives the following algorithm for solving Inverse Problem~1.

\medskip
{\bf Algorithm 1. }{\it Let the spectra $\{\rho_{n,j}\}_{n\in{\mathbb Z}_j}$ of some problems ${\cal
L}_j(q_0,q_1),\,j=0,1,$ be given.

(i) Calculate $\alpha_0,$ $\alpha_1$ and $\omega$ by the formulae (\ref{4.34.1}) and (\ref{4.34.5}), in which
the sequences $\{m_{k,l}\},$ $l=\overline{1,3},$ are chosen so that (\ref{4.34.4}) is fulfilled;

(ii) Construct the functions $w_{j,\nu}(x)\in L_2(0,\pi-a_\nu),$ $j,\nu=0,1,$ in representations (\ref{4.11})
and~(\ref{4.12}) by inverting the corresponding Fourier transforms:
$$
\left[\begin{array}{c}w_{j,j}(x)\\ w_{j,1-j}(x)\end{array}\right] =\frac1\pi\sum_{n=-\infty}^\infty
\theta_j(n) \left[\begin{array}{c}\cos nx\\ \sin nx\end{array}\right], \quad j=0,1,
$$
where the functions $\theta_0(\rho)$ and $\theta_1(\rho)$ are determined by formulae~(\ref{4.33})
and~(\ref{4.34}), respectively, with $\Delta_0(\rho)$ and $\Delta_1(\rho)$ constructed by (\ref{4.13});

(iii) Find the functions $q_0(x)\in L(a_0,\pi)$ and $p(x)\in L(a_1,\pi)$ by formulae~(\ref{4.30.2})
and~(\ref{4.31}) and by formula (\ref{4.17}), respectively, with $w_{j,\nu}(x),\,j,\nu=0,1,$ constructed on
step~(ii);

(iv) Finally, construct the function $q_1(x)\in W_2^1[a_1,\pi]$ by formula (\ref{4.39.1}) or by formula
(\ref{4.34.7}).}

\medskip
{\bf Remark 2.} As in \cite{ButYur19}, step (ii) of Algorithm~1 can be refined by changing to recovering the
functions $w_{j,\nu}(x),\,j,\nu=0,1,$ from certain subspectra depending on the values $a_0$ and~$a_1.$
\\

{\large\bf Appendix A}
\\

Here we obtain an analog of Theorem~4 for the boundary value problems ${\cal L}_j(q),\,j=0,1,$ that consist
of (\ref{0}) and (\ref{2}):
$$
-y''(x)+q(x)y(x-a)=\lambda y(x),\quad 0<x<\pi, \quad y(0)=y^{(j)}(\pi)=0,
$$
where $q(x)\in L_2(0,\pi)$ is a complex-valued function, $q(x)=0$ on $(0,a),$ while $a\in[2\pi/5,\pi).$ Let
$\{\lambda_{n,j}\}_{n\ge1}$ be the spectrum of ${\cal L}_j(q).$ Consider the following inverse problem

\medskip
{\bf Inverse Problem A.} Given $\{\lambda_{n,0}\}_{n\ge1}$ and $\{\lambda_{n,1}\}_{n\ge1},$ find the
potential~$q(x).$

\medskip
As was mentioned in Introduction with references to \cite{ButYur19}, Inverse Problem~A is overdetermined. In
\cite{ButYur19} for $a\in[\pi/2,\pi),$ it was established, in particular, that for unique determination of
$q(x)$ by the subspectra $\{\lambda_{n_k,0}\}_{k\ge1}$ and $\{\lambda_{n_k,1}\}_{k\ge1},$ it is necessary and
sufficient that each of the functional systems $\{\cos n_kx\}_{k\ge1}$ and $\{\sin(n_k-1/2)x\}_{k\ge1}$ is
complete in $L_2(0,\pi-a).$ Moreover, the appropriate asymptotics along with Riesz-basisness of these two
systems is sufficient for solvability of Inverse Problem~A. In particular, solely the asymptotics is a
necessary and sufficient condition of the solvability when this Riesz-basisness is patently the case.
Analogous results can be obtained also for $a\in[2\pi/5,\pi/2).$ So far, these results remain sole ones
dealing with the question of solvability of Inverse Problem~A.

Despite the overdetermination of Inverse Problem~A, one can obtain necessary and sufficient conditions for it
solvability given the full spectra as a particular case of Theorem~4.

\medskip
{\bf Theorem A. }{\it Let $a\in[2\pi/5,\pi).$ Then for any sequences of complex numbers
$\{\lambda_{n,0}\}_{n\ge1}$ and $\{\lambda_{n,1}\}_{n\ge1}$ to be the spectra of some boundary value problems
${\cal L}_0(q)$ and ${\cal L}_1(q),$ respectively, it is necessary and sufficient to satisfy the following
two conditions:

(i) For $j=0,1,$ the following asymptorics holds:
$$
\lambda_{n,j}=\Big(n-\frac j2 +\frac{\omega\cos(n-j/2)a}{\pi n}+\frac{\varkappa_n}n\Big)^2, \quad
\omega\in{\mathbb C},
$$
where, as before, one and the same symbol $\{\varkappa_n\}$ denotes different sequences in $l_2;$

(ii) The exponential types of the functions $\theta_0(\rho)$ and $\theta_1(\rho)$ do not exceed $\pi-a,$
where
\begin{equation}\label{A2}
\theta_0(\rho)=\rho^2\Delta_0(\rho) -\rho\sin\rho\pi +\omega\cos\rho(\pi-a), \;\;
\theta_1(\rho)=\rho\Delta_1(\rho) -\rho\cos\rho\pi -\omega\sin\rho(\pi-a),
\end{equation}
\begin{equation}\label{A3}
\Delta_j(\rho)=\pi^{1-j}\prod_{n=1}^\infty\frac{\lambda_{n,j}-\rho^2}{(n-j/2)^2}, \quad j=0,1.
\end{equation}}

As Corollary~2 from Theorem~4, one can obtain the following corollary from Theorem~A.

\medskip
{\bf Corollary A. }{\it Let $a\in[2\pi/5,\pi).$ Then arbitrary sequences of complex numbers
$\{\lambda_{n,0}\}_{n\ge1}$ and $\{\lambda_{n,1}\}_{n\ge1}$ are the spectra of some boundary value problems
${\cal L}_0(q)$ and ${\cal L}_1(q),$ respectively, if and only if their convergence exponents are equal
to~$1/2,$ and the functions $\theta_j(\rho),\,j=0,1,$ determined by (\ref{A2}) and (\ref{A3}) satisfy the
following conditions:
$$
\theta_j(x)\in L_2(-\infty,\infty), \quad |\theta_j(\rho)|\le C\exp((\pi-a)|\rho|), \quad
\theta_j(-\rho)=(-1)^j\theta_j(\rho).
$$}
\\

{\large\bf Appendix B}
\\

Let $\pi/2\le a_1\le a_2<\pi$ and consider the boundary value problem ${\cal B}$ for the equation
\begin{equation}\label{B.1}
-y''(x)+q_1(x)y(x-a_1)+q_2(x)y(x-a_2)=\lambda y(x),\quad 0<x<\pi,
\end{equation}
where $q_\nu(x)=0$ on $(0,a_\nu)$ and $q_\nu(x)\in L_2(a_\nu,\pi),$ along with the two-point boundary
conditions of the general form:
\begin{equation}\label{B.1.0}
U_\nu(y):=h_{0,\nu}y(0)+h_{1,\nu}y'(0) +H_{0,\nu}y(\pi)+H_{1,\nu}y'(\pi) =0, \quad \nu=1,2,
\end{equation}
with arbitrary complex coefficients $h_{j,\nu}$ and $H_{j,\nu}.$ Let $S(x,\lambda)$ and $C(x,\lambda)$ be
solutions of equation (\ref{B.1}) under the initial conditions $S(0,\lambda)=C'(0,\lambda)=0$ and
$S'(0,\lambda)=C(0,\lambda)=1.$ The spectrum ${\rm sp}({\cal B})$ of the problem ${\cal B}$ coincides with
zeros of its characteristic function
$$
\Delta(\lambda):=\left|\begin{array}{cc}U_1(S)&U_1(C)\\U_2(S)&U_2(C)\end{array}\right|
=\left|\begin{array}{cc}h_{1,1} +H_{0,1}\Delta_0(\lambda)+H_{1,1}\Delta_1(\lambda)& h_{0,1}
+H_{0,1}\Theta_0(\lambda)+H_{1,1}\Theta_1(\lambda)\\ h_{1,2}
+H_{0,2}\Delta_0(\lambda)+H_{1,2}\Delta_1(\lambda)& h_{0,2}
+H_{0,2}\Theta_0(\lambda)+H_{1,2}\Theta_1(\lambda)\end{array}\right|,
$$
where $\Delta_j(\lambda)=S^{(j)}(\pi,\lambda)$ is the characteristic function of the boundary value problem
${\cal B}_j:={\cal B}_j(q_0,q_1)$ consisting of equation (\ref{B.1}) and the boundary conditions (\ref{2}),
while $\Theta_j(\lambda)=C^{(j)}(\pi,\lambda)$ is the characteristic function of the problem for equation
(\ref{B.1}) under the boundary conditions
$$
y'(0)=y^{(j)}(\pi)=0.
$$

The following proposition means, actually, that specification of ${\rm sp}({\cal B})$ gives no additional
information on the potentials $q_1(x)$ and $q_2(x)$ after specifying the spectra ${\rm sp}({\cal B}_0)$
and~${\rm sp}({\cal B}_1).$

\medskip
{\bf Proposition B. }{\it Specification of the functions $\Delta_0(\lambda)$ and $\Delta_1(\lambda)$ along
with the coefficients $h_{j,\nu}$ and $H_{j,\nu}$ for $j=0,1$ and $\nu=1,2$ uniquely determines the function
$\Delta(\lambda).$}

\medskip
{\it Proof.} It is sufficient to show that specification of the functions $\Delta_0(\lambda)$ and
$\Delta_1(\lambda)$ uniquely determines the functions $\Theta_0(\lambda)$ and $\Theta_1(\lambda).$ Indeed,
analogously to Lemma~2 one can obtain the following representations
\begin{equation}\label{B.1.1}
\Delta_0(\lambda)=\Delta_0^0(\lambda) +\sum_{\nu=1}^2\int\limits_0^{\pi-a_\nu} w_{0,\nu}(x) \frac{\cos\rho
x}{\rho^2}\,dx, \quad \Delta_0^0(\lambda):=\frac{\sin\rho\pi}\rho
-\sum_{\nu=1}^2\omega_\nu\frac{\cos\rho(\pi-a_\nu)}{\rho^2},
\end{equation}
\begin{equation}\label{B.1.2}
\Delta_1(\lambda)=\Delta_1^0(\lambda)+\sum_{\nu=1}^2\int\limits_0^{\pi-a_\nu} w_{1,\nu}(x) \frac{\sin\rho
x}\rho\,dx, \quad \Delta_1^0(\lambda):=\cos\rho\pi+\sum_{\nu=1}^2\omega_\nu\frac{\sin\rho(\pi-a_\nu)}\rho,
\end{equation}
where $\rho^2=\lambda$ and
$$
\omega_\nu=\frac12\int\limits_{a_\nu}^\pi q_\nu(x)\,dx, \quad
w_{j,\nu}(x)=\frac14\Big(q_\nu\Big(\frac{\pi+x+a_\nu}2\Big) +(-1)^j q_\nu\Big(\frac{\pi-x+a_\nu}2\Big)\Big).
$$
Moreover, in a similar way, one can get also the representations
$$
\Theta_0(\lambda)=\Delta_1^0(\lambda) -\sum_{\nu=1}^2\int\limits_0^{\pi-a_\nu} w_{1,\nu}(x) \frac{\sin\rho
x}\rho\,dx, \;\; \Theta_1(\lambda)=-\lambda\Delta_0^0(\lambda) +\sum_{\nu=1}^2\int\limits_0^{\pi-a_\nu}
w_{0,\nu}(x) \cos\rho x\,dx.
$$
Thus, we arrive at the relations
$$
\Theta_0(\lambda)=2\Delta_1^0(\lambda)-\Delta_1(\lambda), \quad
\Theta_1(\lambda)=\lambda(\Delta_0(\lambda)-2\Delta_0^0(\lambda)),
$$
which finish the proof because $\Delta_j^0(\lambda)$ is determined by specifying $\Delta_j(\lambda).$
$\hfill\Box$

\medskip
Obviously, no set of spectra determines the functions $q_1(x)$ and $q_2(x)$ separately if $a_1=a_2.$ The
following example shows that they cannot be completely distinguished also when $a_1\ne a_2.$

\medskip
{\bf Example B.} Let
\begin{equation}\label{B.2}
q_1(x)=\left\{\begin{array}{rc}0, & \displaystyle
a_1<x<\frac{a_1+a_2}2,\\[3mm]
1, & \displaystyle \frac{a_1+a_2}2<x<\frac{a_1+\pi}2,\\[3mm]
-1, & \displaystyle \frac{a_1+\pi}2<x<\pi-\frac{a_2-a_1}2,\\[3mm]
0, & \displaystyle \pi-\frac{a_2-a_1}2<x<\pi,
\end{array}\right. \quad
q_2(x)=\left\{\begin{array}{rc}
-1, & \displaystyle a_2<x<\frac{a_2+\pi}2,\\[3mm]
1, & \displaystyle \frac{a_2+\pi}2<x<\pi,
\end{array}\right.
\end{equation}
and, hence,
\begin{equation}\label{B.3}
q_2(x)=-q_1\Big(x-\frac{a_2-a_1}2\Big), \quad a_2<x<\pi.
\end{equation}
Then the spectra of the problems ${\cal B}_0(q_0,q_1)$ and ${\cal B}_1(q_0,q_1)$ coincide with the ones of
${\cal B}_0(0,0)$ and ${\cal B}_1(0,0),$ respectively. Indeed, according to the relations
$$
{\mathfrak L}(\rho):=\Delta_1(\lambda)+i\rho\Delta_0(\lambda), \quad \Delta_0(\lambda)=\frac{{\mathfrak
L}(\rho)-{\mathfrak L}(-\rho)}{2i\rho}, \quad \Delta_1(\lambda)=\frac{{\mathfrak L}(\rho)+{\mathfrak
L}(-\rho)}2,
$$
specification of both spectra is equivalent to specification of the function ${\mathfrak L}(\rho),$ which, in
turn, is the characteristic function of the Regge-type problem for equation (\ref{B.1}) along with the
boundary conditions
\begin{equation}\label{B.4}
y(0)=y'(\pi)+i\rho y(\pi)=0.
\end{equation}
On the other hand, the following representation holds (see, e.g., \cite{Sh19}):
\begin{equation}\label{B.5}
{\mathfrak L}(\rho)\exp(-i\rho\pi)-1=\sum_{\nu=1}^2\frac{\omega_\nu}{i\rho}\exp(-i\rho a_\nu)
-\sum_{\nu=1}^2\frac{\exp(i\rho a_\nu)}{2i\rho}\int\limits_{a_\nu}^\pi q_\nu(x)\exp(-2i\rho x)\,dx,
\end{equation}
which also can be easily obtained by using (\ref{B.1.1}) and (\ref{B.1.2}). According to (\ref{B.2}), we have
$$
2\omega_1=\int\limits_{\frac{a_1+a_2}2}^{\frac{a_1+\pi}2}dx
-\int\limits_{\frac{a_1+\pi}2}^{\pi-\frac{a_2-a_1}2}dx=0, \quad
2\omega_2=-\int\limits_{a_2}^\frac{a_2+\pi}2dx +\int\limits_\frac{a_2+\pi}2^\pi dx=0.
$$
Moreover, by virtue of (\ref{B.2}) and (\ref{B.3}), we get
$$
\sum_{\nu=1}^2\exp(i\rho a_\nu)\int\limits_{a_\nu}^\pi q_\nu(x)\exp(-2i\rho x)\,dx  =\exp(i\rho
a_1)\int\limits_{\frac{a_1+a_2}2}^{\pi-\frac{a_2-a_1}2} q_1(x)\exp(-2i\rho x)\,dx \qquad\qquad\quad
$$
$$
\qquad\qquad\qquad\qquad\qquad\qquad\qquad\qquad -\exp(i\rho a_2)\int\limits_{a_2}^\pi
q_1\Big(x-\frac{a_2-a_1}2\Big)\exp(-2i\rho x)\,dx =0.
$$
Hence, according to (\ref{B.5}), we have ${\mathfrak L}(\rho)=\exp(i\rho\pi),$ which is the characteristic
function of the problem (\ref{B.1}), (\ref{B.4}) with the zero potentials. Thus, specification of the spectra
of the problems ${\cal B}_0(q_1,q_2)$ and ${\cal B}_1(q_1,q_2)$ does not uniquely determine the functions
~$q_1(x)$ and~$q_2(x).$

\medskip
Finally, we note that Example~B refutes both Theorem~3.1 in \cite{Sh19} and Theorem~3.1 in \cite{Sh20}.
Moreover, according to this counterexample along with Proposition~B, the functions $q_1(x)$ and $q_2(x)$
cannot be uniquely determined by specifying any set of the spectra of boundary value problems having the form
(\ref{B.1}) and (\ref{B.1.0}).
\\

{\bf Acknowledgement.} This research was partially supported by the Ministry of Science and Technology of
Taiwan under Grant no. 108-2115-M-032-005-. The first author was supported by Grant 20-31-70005 of the
Russian Foundation for Basic Research. The second authors was supported by Grant 1.1660.2017/4.6 of the
Ministry of Science and Higher Education of the Russian Federation.


\begin{thebibliography}{99}

\bibitem{Pik91}
Pikula M. {\it Determination of a Sturm--Liouville-type differential operator with delay argument from two
spectra}, Mat. Vestnik 43 (1991) no.3-4, 159--171.

\bibitem{FrYur12}
Freiling G. and Yurko V.A. {\it Inverse problems for Sturm--Liouville differential operators with a constant
delay}, Appl. Math. Lett. 25 (2012) 1999--2004.

\bibitem{Yang14}
Yang C.-F. {\it Inverse nodal problems for the Sturm--Liouville operator with a constant delay}, J. Diff.
Eqns. 257 (2014) no.4, 1288--1306.

\bibitem{VladPik16}
Vladi\v{c}i\'c V. and Pikula M. {\it An inverse problem for Sturm--Liouville-type differential equation with
a constant delay}, Sarajevo J. Math. 12 (2016) no.1, 83--88.

\bibitem{ButYur19}
Buterin S.A. and Yurko V.A. {\it An inverse spectral problem for Sturm--Liouville operators with a large
delay}, Anal. Math. Phys. 9 (2019) no.1, 17--27. (published online in 2017)

\bibitem{ButPikYur17}
Buterin S.A., Pikula M. and Yurko V.A. {\it Sturm--Liouville differential operators with deviating argument},
Tamkang J. Math. 48 (2017) no.1, 61--71.

\bibitem{Ign18}
Ignatiev M.Yu. {\it On an inverse Regge problem for the Sturm--Liouville operator with deviating argument},
J. Samara State Tech. Univ., Ser. Phys. Math. Sci. 22 (2018) no.2, 203--211.

\bibitem{BondYur18-1}
Bondarenko N. and Yurko V. {\it An inverse problem for Sturm--Liouville differential operators with deviating
argument}, Appl. Math. Lett. 83 (2018) 140--144.

\bibitem{BondYur18-2}
Bondarenko N.P. and Yurko V.A. {\it Partial inverse problems for the Sturm--Liouville equation with deviating
argument}, Math. Meth. Appl. Sci. 41 (2018) 8350--8354.

\bibitem{VPV19}
Vladi\v{c}i\'c V., Pikula M. and Vojvodi\'c B. {\it Inverse spectral problems for Sturm--Liouville operators
with a constant delay less than half the length of the interval and Robin boundary conditions}, Results Math.
(2019) 74:45.

\bibitem{DV19}
Djuri\'c N. and Vladi\v{c}i\'c V. {\it Incomplete inverse problem for Sturm--Liouville type differential
equation with constant delay}, Results Math. (2019) 74:161.

\bibitem{SS19}
Sat M. and Shieh C.-T. {\it Inverse nodal problems for integro-differential operators with a constant delay},
J. Inverse Ill-Posed Probl. 27 (2019) no.4, 501--509.

\bibitem{Yur19}
Yurko V.A. {\it An inverse spectral problem for second order differential operators with retarded argument},
Results Math. (2019) 74:71.

\bibitem{Yur19-2}
Yurko V.A. {\it Recovering differential operators with a retarded argument}, Diff. Eqns. 55 (2019) no.4,
524--528.

\bibitem{WShM19}
Wang Y.P., Shieh C.T. and Miao H.Y. {\it Reconstruction for Sturm--Liouville equations with a constant delay
with twin-dense nodal subsets}, Inv. Probl. Sci. Eng. 27 (2019) no.5, 608--617.

\bibitem{Yur20}
Yurko V.A. {\it Solution of Inverse Problems for Differential Operators with Delay}, in V.~Kravchenko and
S.~Sitnik (Eds.), Trends in Mathematics: Transmutation Operators and Applications, Birkh\"auser, Basel, 2020.
P.~467--475.

\bibitem{DurBut}
Djuri\'c N. and Buterin S. {\it On an open question in recovering Sturm--Liouville-type operators with
delay}, Appl. Math. Lett. (2021) 106862 (arXiv:2009.02636, 6pp.)

\bibitem{PVP16}
Vojvodi\'c B., Pikula M. {\it The boundary value problem for the differential of operator Sturm--Liouville
type with $N$ constant delays and asymptotics of eigenvalues}, Math. Montisnigri 35 (2016) 5--21.

\bibitem{Sh19}
Shahriari M. {\it Inverse problem for Sturm--Liouville differential operators with two constant delays},
Turkish J. Math. 43 (2019) 965--976.

\bibitem{Sh20}
Shahriari M. {\it Inverse problem for Sturm--Liouville differential operators with finite number of constant
delays}, Turkish J. Math. 44 (2020) 778--790.

\bibitem{VPV20}
Vojvodi\'c B., Vladi\v{c}i\'c V. and Pikula M. {\it Inverse problems for Sturm--Liouville differential
operators with two constant delays under Robin boundary conditions}, Results Appl. Math.~5 (2020) 100082.

\bibitem{VPVC20}
Vojvodi\'c B., Vladi\v{c}i\'c V., Pikula M. and \c{C}etinkaya F.A. {\it Inverse problems for differential
operators with two delays larger than half the length of the interval and Dirichlet conditions}, Turkish J.
Math. 44 (2020) no.3, 900--905.

\bibitem{BM19}
Buterin S.A. and Malyugina M.A. {\it On recovering differential pencils with delay}, in Sovr. Met. Teor.
Funk. i Smezhn. Vopr., Voronesh Univ. Press, Voronezh, 2019, 64--66.

\bibitem{BMS19}
Buterin S.A., Malyugina M.A. and Shieh C.-T. {\it An inverse spectral problem for functional-differential
pencils with delay}, in Sovr. Probl. Mat. i Mekh., MAKS Press, Moscow, 2019, 422--425.

\bibitem{ButYur12}
Buterin S.A. and Yurko V.A. {\it Inverse problems for second-order differential pencils with Dirichlet
boundary conditions}, J. Inverse Ill-Posed Probl. 20 (2012) 855--881.

\bibitem{But20}
Buterin S. {\it On a Transformation Operator Approach in the Inverse Spectral Theory of Integral and
Integro-Differential Operators}, in V.~Kravchenko and S.~Sitnik (Eds.), Trends in Mathematics: Transmutation
Operators and Applications, Birkh\"auser, Basel, 2020. P.~337--367.



\bibitem{FY01} Freiling G. and Yurko V.A. {\it Inverse Sturm--Liouville Problems and Their Applications},
NOVA Science Publishers, New York, 2001.

\bibitem{BFY14}
Buterin S.A., Freiling G., Yurko V.A. {\it Lectures in the Theory of Entire Functions}, Schriftenriehe der
Fakult\"at f\"ur Mathematik, Universit\"at Duisbug--Essen, SM-UDE-779, 2014.
\end{thebibliography}
\end{document}